\documentclass[a4paper,reqno,10pt]{amsart}

\usepackage{amsmath,amssymb,amsfonts,amsthm,mathrsfs}
\usepackage{verbatim,wasysym,cite}

\usepackage[latin1]{inputenc}
\usepackage{microtype}
\usepackage{color,enumitem,graphicx}

\usepackage[
colorlinks=true,
urlcolor=blue,
citecolor=red,
linkcolor=blue,
linktocpage,
pdfpagelabels,
bookmarksnumbered,
bookmarksopen
]{hyperref}

\usepackage[english]{babel}
\usepackage[symbol]{footmisc}

\renewcommand{\epsilon}{\varepsilon}

\numberwithin{equation}{section}

\newtheorem{theorem}{Theorem}[section]
\newtheorem{lemma}[theorem]{Lemma}
\newtheorem{remark}[theorem]{Remark}
\newtheorem{definition}[theorem]{Definition}
\newtheorem{proposition}[theorem]{Proposition}

\newcommand{\C}{\mathbb C}

\newcommand{\R}{\mathbb R}

\def\({\left(}
\def\){\right)}
\def\<{\left\langle}
\def\>{\right\rangle}

\DeclareMathOperator{\RE}{Re}

\begin{document}

\title[Inhomogeneous logarithmic Schr\"odinger equation]
{Orbital stability and instability of standing waves
for an inhomogeneous logarithmic Schr\"odinger equation}

 \author[Alex H. Ardila]{Alex H. Ardila}
\address{Alex H. Ardila
		\newline \indent Department of Mathematics, Universidad del Valle, Colombia} 
\email{ardila@impa.br}

\begin{abstract}
{We study} the inhomogeneous logarithmic
Schr\"odinger equation
\[
i\partial_tu+\Delta u+|x|^{-1}u\log|u|^2=0,
\qquad (t,x)\in\R\times\R^N,\quad N\ge2.
\]
{For every initial datum in $H^1(\R^N)$,
we establish global existence of weak solutions.
We also construct explicit positive stationary
solutions and solve the associated Nehari
minimization problem.}
For $\omega\ge0$, the minimizer is unique up to
{a constant} phase and is given by
\[
e^{\frac{N-1}{2}(1+\sqrt{1+\omega})}
e^{-(1+\sqrt{1+\omega})|x|},
\]
whereas for $\omega<0$ the infimum is zero and is not
attained. Finally, we prove that the standing waves
associated with these profiles are orbitally stable
for $\omega>0$ by a variational argument, and orbitally
unstable for $-1<\omega<-1+(N-1)^{-2}$
by combining spectral analysis,
{localized energy and mass estimates,
and phase modulation}.
\end{abstract}


\maketitle

\medskip

\section{Introduction}
\label{intr}
In this paper, we study the inhomogeneous logarithmic
Schr\"odinger equation
\begin{equation}
i\partial_tu+\Delta u+|x|^{-1}u\log|u|^2=0,
\qquad (t,x)\in\R\times\R^N,\quad N\ge2,
\label{NLS}
\end{equation}
where $u:\R\times\R^N\to\C$. Our aim is to characterize the
minimizers of the associated Nehari problem and to study
the orbital stability and instability of  standing
waves. 
Logarithmic Schr\"odinger equations with inhomogeneous
coefficients have been {studied in
connection with quantum liquids and localized wave
propagation}; see \cite{Z18,Z19,CABC} and the references therein.
Zloshchastiev~\cite{Z18,Z19} relates the logarithmic
coupling to temperature and
{examines how temperature gradients
affect quantum liquids}.
Cala\c{c}a, Avelar, Bazeia, and Cardoso~\cite{CABC}
{study how variable linear and nonlinear
coefficients can be used to modulate localized solutions}.

The analysis of logarithmic Schr\"odinger equations goes
back to Cazenave and Haraux~\cite{CH80}, who established
existence and uniqueness results in suitable functional
spaces. For the classical logarithmic Schr\"odinger equation,
Cazenave~\cite{C83} studied the stability of Gaussian
standing waves under radial perturbations, and
Ardila~\cite{A16} proved orbital stability in the natural
energy space without a radial restriction. Variational
questions for the stationary equation were studied by
d'Avenia, Montefusco, and Squassina~\cite{DMS}, who used
nonsmooth critical point theory to obtain infinitely
many weak solutions and proved uniqueness, up to
translations, and nondegeneracy of the positive solution.

The influence of defects, geometry, and external potentials
has been investigated in several settings. 
See, for example, \cite{AA18,A17,AG17,G19, CS22} and the references therein.
In particular, we mention the work of Carles and Su~\cite{CS22}, who showed that a repulsive harmonic potential can produce distinct positive Gaussian stationary profiles whose standing waves are orbitally unstable.

A related class of problems is given by the inhomogeneous
nonlinear Schr\"odinger equation
\[
i\partial_tu+\Delta u+\lambda |x|^{-b}|u|^\alpha u=0,
\qquad \lambda=\pm1,
\]
whose standing waves have been extensively studied in
the focusing case $\lambda=1$.
For more general coefficients behaving like $|x|^{-b}$
at infinity, de Bouard and Fukuizumi~\cite{BF05}
established orbital stability at small positive frequencies
in the mass-subcritical regime, whereas Fukuizumi and
Ohta~\cite{FO05} proved instability in the
mass-supercritical, energy-subcritical regime, under
suitable hypotheses on the coefficient.
Genoud and Stuart~\cite{GS08} constructed branches of
standing waves and obtained stability and instability
results at low frequencies.
For the coefficient $|x|^{-b}$, with $0<b<\min\{2,N\}$,
Genoud~\cite{Ge12} proved strong instability by blow-up
at the mass-critical exponent $\alpha=(4-2b)/N$.
Ardila and Cardoso~\cite{AC21} established orbital
stability for $0<\alpha<(4-2b)/N$ and strong instability
by blow-up in the mass-supercritical, energy-subcritical
range using variational methods.
Ramadan and Stefanov~\cite{RS20} obtained existence,
nondegeneracy, and stability results in a framework
including classical and fractional dispersion, under
additional restrictions on the parameters.

The energy and mass associated with \eqref{NLS}
are defined, respectively, by
\begin{equation*}
E(u)
=
\frac12\int_{\R^N}|\nabla u|^2\,dx
-\frac12\int_{\R^N}
\frac{|u|^2\log|u|^2-|u|^2}{|x|}\,dx,
\qquad
M(u)=\int_{\R^N}|u|^2\,dx.
\end{equation*}
By Lemma~\ref{le01} below, $E$ is
of class $C^1$ on $H^1(\R^N)$, and
\[
E'(u)=-\Delta u-|x|^{-1}u\log|u|^2
\quad\text{in }H^{-1}(\R^N).
\]
Moreover, $E'$ is continuous and maps bounded sets
into bounded sets. Equation~\eqref{NLS} can therefore
be interpreted in $H^{-1}(\R^N)$ at the regularity
stated below. Although mass and energy are formally
conserved, our construction of weak solutions gives
mass conservation and an energy inequality. These
properties suffice for the stability/instability argument.

\begin{proposition}
\label{pr01}
Let $N\ge2$. For every $u_0\in H^1(\R^N)$, there
exists a global weak solution $u\in L^\infty(\R;H^1(\R^N))
\cap C_w(\R;H^1(\R^N))
\cap W_{\mathrm{loc}}^{1,\infty}
(\R;H^{-1}(\R^N))$
of \eqref{NLS} such that $u(0)=u_0$ and $ \sup_{t\in\R}\|u(t)\|_{H^1}
<\infty$. Furthermore, conservation of mass
and the energy inequality hold: $M(u(t))=M(u_0)$, $E(u(t))\le E(u_0)$ for all $t\in\R$.
In addition, $\lim_{t\to0}\|u(t)-u_0\|_{H^1}=0$.
\end{proposition}

Proposition~\ref{pr01} provides global weak solutions
for arbitrary initial data in $H^1(\R^N)$.
The following result establishes uniqueness whenever
one of the solutions satisfies an exponential lower
bound on each finite time interval. This condition
is satisfied, in particular, by the explicit standing
waves considered below (cf. Remark \ref{re01}).
\begin{proposition}
\label{pr04}
Let $N\ge2$, $T>0$, and $u_0\in H^1(\R^N)$.
Let $u$ and $v$ be two solutions of \eqref{NLS}
with initial datum $u_0$, satisfying the properties
stated in Proposition~\ref{pr01}.
Assume that there exists $A_T>0$ such that
\begin{equation}
|v(t,x)|\ge e^{-A_T(1+|x|)}
\quad\text{for almost every }
(t,x)\in(0,T)\times\R^N.
\label{un01}
\end{equation}
Then
\begin{equation}
u(t)=v(t)\quad\text{in }H^1(\R^N)
\quad\text{for every }t\in[0,T].
\label{un02}
\end{equation}
\end{proposition}

We now turn to standing waves of the form
$u(t,x)=e^{i\omega t}\phi(x)$, where $\omega\in\R$.
Their profiles solve the stationary problem
\begin{equation}
\begin{cases}
-\Delta\phi+\omega\phi-|x|^{-1}\phi\log|\phi|^2=0
& \text{in }H^{-1}(\R^N),\\
\phi\in H^1(\R^N)\setminus\{0\}.
\end{cases}
\label{GS}
\end{equation}
To formulate the variational problem, we introduce
the action and the Nehari functional.

\begin{definition}
\label{de01}
For $\omega\in\R$, we define the following functionals
of class $C^1$ on $H^1(\R^N)$:
\[
\begin{aligned}
S_\omega(u)
&:=
\frac12\left\|\nabla u\right\|_{L^2(\R^N)}^2
+\frac{\omega}{2}\left\|u\right\|_{L^2(\R^N)}^2
-\frac12\int_{\R^N}
\frac{|u|^2\log|u|^2-|u|^2}{|x|}\,dx,
\end{aligned}
\]
and
\[
I_\omega(u)
:=
\left\|\nabla u\right\|_{L^2(\R^N)}^2+\omega\left\|u\right\|_{L^2(\R^N)}^2
-\int_{\R^N}\frac{|u|^2\log|u|^2}{|x|}\,dx.
\]
The functional $S_\omega$ is called the action, and
$I_\omega(u)=\langle S_\omega'(u),u\rangle$ is the Nehari functional.
\end{definition}

Notice that the stationary equation in \eqref{GS}
is equivalent to $S_\omega'(\phi)=0$. A nonzero
solution $\phi\in H^1(\R^N)$ of \eqref{GS} is called
a ground state if it has minimal action among all
nonzero solutions of \eqref{GS}. To study this
variational structure, we consider the minimization
problem
\begin{equation}
\begin{aligned}
d(\omega)
&:=
\inf\left\{
S_\omega(u):
u\in H^1(\R^N)\setminus\{0\},\
I_\omega(u)=0
\right\}\\
&=
\tfrac12\inf\left\{
\int_{\R^N}\frac{|u|^2}{|x|}\,dx:
u\in H^1(\R^N)\setminus\{0\},\
I_\omega(u)=0
\right\},
\end{aligned}
\label{MPE}
\end{equation}
and denote its set of minimizers by
\[
\mathcal N_\omega
:=
\left\{
\phi\in H^1(\R^N)\setminus\{0\}:
S_\omega(\phi)=d(\omega),\
I_\omega(\phi)=0
\right\}.
\]
The set $\left\{
u\in H^1(\R^N)\setminus\{0\}:
I_\omega(u)=0
\right\}$
is called the Nehari manifold and contains every
nonzero critical point of $S_\omega$.
Whenever $\mathcal N_\omega$ is nonempty, it coincides
with the set of ground states. Indeed, on the Nehari
manifold,
\[
\langle I_\omega'(u),u\rangle
=
-2\int_{\R^N}\frac{|u|^2}{|x|}\,dx<0.
\]
The Lagrange multiplier rule therefore shows that
every minimizer of \eqref{MPE} is a critical point
of $S_\omega$, and it has minimal action among all
nonzero stationary solutions. Conversely, when the
infimum is attained, every ground state attains
$d(\omega)$.

\begin{proposition}
\label{pr02}
Let $N\ge2$. The stationary problem \eqref{GS} admits
the following positive solutions:
\begin{enumerate}[label=\textup{(\roman*)}]
\item If $\omega=-1$, then $\phi(x)=e^{\frac{N-1}{2}}e^{-|x|}.$

\item If $-1<\omega<0$, then $\phi(x)
=
e^{\frac{N-1}{2}(1+\sqrt{1+\omega})}
e^{-(1+\sqrt{1+\omega})|x|}$
and $\psi(x)
=
e^{\frac{N-1}{2}(1-\sqrt{1+\omega})}
e^{-(1-\sqrt{1+\omega})|x|}$
are two distinct solutions.

\item If $\omega\ge0$, then $\phi(x)
=
e^{\frac{N-1}{2}(1+\sqrt{1+\omega})}
e^{-(1+\sqrt{1+\omega})|x|}.$
\end{enumerate}
\end{proposition}

\begin{remark}
\label{re01}
Let $\phi$ be any of the profiles in
Proposition~\ref{pr02}, with corresponding frequency
$\omega$. Then $u(t,x)=e^{i\omega t}\phi(x)$
is the unique solution in the class of Proposition~\ref{pr01}
with initial
datum $\phi$. Indeed, the explicit formulas imply that, for every
$T>0$, there exists $A_T>0$ such that
\[
|u(t,x)|=\phi(x)\ge e^{-A_T(1+|x|)}
\qquad\text{for all }(t,x)\in[0,T]\times\R^N.
\]
Proposition~\ref{pr04} therefore gives uniqueness
on $[0,T]$. Since $T>0$ is arbitrary, uniqueness
holds for all times, using time reversal for negative times.
\end{remark}

\begin{theorem}
\label{th01}
Let $N\ge2$ and $\omega\in\R$. Then the following
assertions hold.
\begin{enumerate}[label=\textup{(\roman*)}]
\item If $\omega\ge0$, the infimum in \eqref{MPE}
is attained, and
\[
\mathcal N_\omega
=
\left\{
e^{i\theta}
e^{\frac{N-1}{2}(1+\sqrt{1+\omega})}
e^{-(1+\sqrt{1+\omega})|x|}:
\theta\in\R
\right\}.
\]
Moreover,
\[
d(\omega)
=
\frac{|\mathbb S^{N-1}|\Gamma(N-1)}
{2[2(1+\sqrt{1+\omega})]^{N-1}}
e^{(N-1)(1+\sqrt{1+\omega})},
\]
where $\Gamma$ denotes the gamma function,
and $|\mathbb S^{N-1}|$ is the surface measure
of the unit sphere in $\R^N$.
\item If $\omega<0$, then $d(\omega)=0$ and the minimization problem \eqref{MPE} has no minimizers.

\end{enumerate}
\end{theorem}

Equation \eqref{NLS} is invariant under phase rotations,
whereas the spatial weight $|x|^{-1}$ breaks translation
invariance. Accordingly, orbital stability is defined
with respect to the phase orbit. Since
Proposition~\ref{pr01} does not assert uniqueness for
arbitrary initial data, we require the stability
property to hold for every solution in the class
specified there.

\begin{definition}
\label{de02}
A standing-wave solution $e^{i\omega t}\phi(x)$ of
\eqref{NLS} is said to be orbitally stable in
$H^1(\R^N)$ if, for every $\varepsilon>0$, there
exists $\delta>0$ such that, whenever
$u_0\in H^1(\R^N)$ satisfies $\|u_0-\phi\|_{H^1}<\delta$,
every global weak solution $u$ with initial datum
$u_0$ satisfying the properties stated in
Proposition~\ref{pr01} obeys
\[
\sup_{t\in\R}\inf_{\theta\in\R}
\|u(t)-e^{i\theta}\phi\|_{H^1}<\varepsilon.
\]
Otherwise, the standing wave is said to be orbitally
unstable in $H^1(\R^N)$.
\end{definition}

Although the energy functional $E$ belongs to
$C^1$, it is not of class $C^2$.
Consequently, the abstract stability theory of
Grillakis, Shatah, and Strauss~\cite{GSS1}
cannot be applied directly in this space.
We therefore establish orbital stability through
a variational argument based on the Nehari
minimization problem and compactness. For orbital
instability, we use an ejection argument: spectral
coercivity, phase modulation, and localized estimates
control the nonlinear errors and force the unstable
coordinate to grow until the solution leaves a
fixed neighborhood of the standing-wave orbit.

Our next result establishes the orbital stability
of the ground states for every positive frequency.

\begin{theorem}
\label{th02}
Let $N\ge2$ and $\omega>0$. Then the standing wave $u(t,x)
=
e^{i\omega t}
e^{\frac{N-1}{2}(1+\sqrt{1+\omega})}
e^{-(1+\sqrt{1+\omega})|x|}$
is orbitally stable in $H^1(\R^N)$.
\end{theorem}

The proof of Theorem~\ref{th02} combines the
variational characterization of the ground states
in Theorem~\ref{th01} with the compact weighted
embeddings of Lemma~\ref{CCD}. These embeddings
yield compactness of minimizing sequences up to
a constant phase. The strict convexity of
$d(\omega)$ for $\omega>0$, together with mass
conservation and the energy inequality from
Proposition~\ref{pr01}, controls the weighted mass $\int_{\R^N}|x|^{-1}|u|^2\,dx$
along solutions starting near the standing-wave
orbit. If stability failed, these estimates would
produce a sequence of solution profiles approaching
the Nehari minimum while remaining separated from
the orbit, contradicting compactness.

We now present our instability result for negative frequencies.

\begin{theorem}
\label{thm1}
Let $N\ge2$ and $-1<\omega<-1+\frac{1}{(N-1)^2}$.
Then the standing wave $u(t,x)
=
e^{i\omega t}
e^{\frac{N-1}{2}(1+\sqrt{1+\omega})}
e^{-(1+\sqrt{1+\omega})|x|}$
is orbitally unstable in $H^1(\R^N)$.
\end{theorem}

The proof of Theorem~\ref{thm1} uses the growth of
a linearly unstable mode to establish orbital
instability, a principle studied for nonlinear
Schr\"odinger equations by Georgiev and
Ohta~\cite{GO12}. We implement this principle through
an ejection argument in the spirit of
Nakanishi and Schlag~\cite[Lemma~3.3]{NS12}
and Nakanishi and Roy~\cite[Proposition~3.7]{NR16}.
We construct arbitrarily small radial perturbations
with the same mass as the standing wave and strictly
smaller energy. Spectral coercivity and phase
modulation reduce the analysis to the stable and
unstable coordinates and a localized mass term.
Localized energy estimates and mass identities
control the nonlinear errors associated with the
logarithmic nonlinearity and the singular coefficient.
A continuity argument then shows that the unstable
coordinate grows exponentially for as long as the
solution remains close to the standing-wave orbit.
This growth forces departure from a fixed
neighborhood and establishes orbital instability
under radial perturbations in $H^1(\R^N)$.

\begin{remark}\label{re02}
The following cases are not covered by our orbital
stability and instability results:
\begin{enumerate}[label=\textup{(\roman*)}]
\item
The upper-branch standing waves
\[
u(t,x)=
e^{i\omega t}
e^{\frac{N-1}{2}(1+\sqrt{1+\omega})}
e^{-(1+\sqrt{1+\omega})|x|},
\qquad
-1+{(N-1)^{-2}}\le\omega\le0.
\]

\item
The lower-branch standing waves
\[
u(t,x)=
e^{i\omega t}
e^{\frac{N-1}{2}(1-\sqrt{1+\omega})}
e^{-(1-\sqrt{1+\omega})|x|},
\qquad -1<\omega<0.
\]

\item
The standing wave at $\omega=-1$, where the two
branches coincide:
\[
u(t,x)=e^{-it}e^{\frac{N-1}{2}}e^{-|x|}.
\]
\end{enumerate}

When $N=2$, the interval in \textup{(i)} reduces
to the single frequency $\omega=0$.
Thus Theorems~\ref{th02} and~\ref{thm1} establish
orbital stability of the upper branch for
$\omega>0$ and orbital instability for
$-1<\omega<0$, leaving the threshold solution
$u(t,x)=e^{1-2|x|}$ unresolved.
For $N\ge2$, the profile at $\omega=0$
is a Nehari minimizer by Theorem~\ref{th01},
but its orbital stability is not established
by Theorem~\ref{th02}, which requires $\omega>0$.
\end{remark}

The paper is organized as follows.
Section~\ref{prel} establishes the $C^1$ regularity
of the energy functional, the compact weighted
embeddings, and the sharp weighted logarithmic
Sobolev inequality with its equality cases.
For $N\ge3$, the latter follows from the results
of Dolbeault and Zuniga~\cite{DZ22}, while the
case $N=2$ is treated directly.
Section~\ref{caup} proves the global existence
and conditional uniqueness results in
Propositions~\ref{pr01} and~\ref{pr04}.
Existence is obtained by regularizing both the
singular coefficient and the logarithmic
nonlinearity, following Cazenave~\cite{C03},
whereas uniqueness under \eqref{un01} follows
from nonlinear estimates and Strichartz
estimates~\cite{KT98}.
Section~\ref{vari} constructs the explicit stationary
profiles and characterizes the Nehari minimizers,
proving Proposition~\ref{pr02} and
Theorem~\ref{th01}.

Section~\ref{stab} proves Theorem~\ref{th02}
by a variational argument following
Fukuizumi, Ohta, and Ozawa~\cite{FOO8}.
The proof combines compactness of minimizing
sequences and strict convexity of $d(\omega)$
with mass conservation and the energy inequality.
Section~\ref{inst} establishes Theorem~\ref{thm1}
through an ejection argument based on phase
modulation and localized energy and mass estimates.
These estimates control the nonlinear errors
and ensure growth of the unstable coordinate
until the solution leaves a fixed neighborhood
of the standing-wave orbit.
The spectral properties, coercivity estimates,
and exponential decay of the spectral modes
used in this argument are proved in
Appendix~\ref{asp}.

\subsection*{Notation}

We use the usual Lebesgue and Sobolev spaces
$L^p(\R^N)$ and $H^1(\R^N)$, with
$\|u\|_{H^1}^2:=\left\|u\right\|_{L^2(\R^N)}^2+\left\|\nabla u\right\|_{L^2(\R^N)}^2$.
We write $H^{-1}$ for the dual of $H^1$
under the $L^2$ pairing, and
$C_w(I;H^1)$ for the space of weakly continuous
$H^1$-valued functions on an interval $I$.
The $L^2$ inner product and the weighted mass are
\[
(u,v)_{L^2}:=\operatorname{Re}
\int_{\R^N}u\overline v\,dx
\qquad
m_1(u):=\int_{\R^N}|x|^{-1}|u|^2\,dx.
\]
For $R>0$, we set $B_R:=\{x\in\R^N:|x|<R\}$,
so that $B_1$ is the unit ball.
We write $\log_+s:=\max\{\log s,0\}$ for $s>0$,
with $\log_+0:=0$, and adopt the convention
$0\log0:=0$.
The letter $C$ denotes a positive constant that
may change from line to line.

\section{Preliminaries}
\label{prel}

In this section, we establish the analytic estimates
needed for the Cauchy problem and the variational
analysis of standing waves. We first prove that the
energy functional is well defined and of class $C^1$
on $H^1(\R^N)$. We then establish a sharp weighted
logarithmic Sobolev inequality and characterize its
equality cases, which will be used to identify the
minimizers of the Nehari problem.

\begin{lemma}
\label{le01}
Let $N\ge2$. The energy functional $E:H^1(\R^N)\to\R$ is well
defined and of class $C^1$. Moreover, its
Fr\'echet derivative is given by
\[
E'(u)=-\Delta u-|x|^{-1}u\log|u|^2
\quad\text{in }H^{-1}(\R^N).
\]
\end{lemma}

\begin{proof}
Fix $0<\varepsilon<2/N$. We first establish the
weighted estimate (cf. \cite[Theorem~1.1(i)]{R12})
\begin{equation}
\int_{\R^N}
\frac{|u|^{2-\varepsilon}+|u|^{2+\varepsilon}}{|x|}\,dx
\le
C_{N,\varepsilon}
\left(
\|u\|_{H^1}^{2-\varepsilon}
+\|u\|_{H^1}^{2+\varepsilon}
\right).
\label{eq03}
\end{equation}
On $B_1=\{x:|x|<1\}$,  H\"older's inequality and the
Sobolev embeddings  give
\[
\int_{B_1}\frac{|u|^{2\pm\varepsilon}}{|x|}\,dx
\le
\left\||x|^{-1}\right\|_{L^r(B_1)}
\left\|u\right\|_{L^{(2\pm\varepsilon)r'}(\R^N)}^{2\pm\varepsilon}
\le
C_{N,\varepsilon}\|u\|_{H^1}^{2\pm\varepsilon},
\] where we choose $1<r<N$ sufficiently close to $N$
so that both $(2-\varepsilon)r'$ and
$(2+\varepsilon)r'$ belong to the Sobolev embedding
range, with $r'=r/(r-1)$. This is possible because
$\varepsilon<2/N$.
On the complement of $B_1$, we have
\[
\int_{|x|>1}\frac{|u|^{2+\varepsilon}}{|x|}\,dx
\le
\left\|u\right\|_{L^{2+\varepsilon}(\R^N)}^{2+\varepsilon}
\le
C_{N,\varepsilon}\|u\|_{H^1}^{2+\varepsilon},
\]
whereas H\"older's inequality yields
\[
\begin{aligned}
\int_{|x|>1}\frac{|u|^{2-\varepsilon}}{|x|}\,dx
&\le
\left\|u\right\|_{L^2(\R^N)}^{2-\varepsilon}
\left(
\int_{|x|>1}|x|^{-2/\varepsilon}\,dx
\right)^{\varepsilon/2}
\le
C_{N,\varepsilon}\left\|u\right\|_{L^2(\R^N)}^{2-\varepsilon}.
\end{aligned}
\]
This proves \eqref{eq03}. The elementary estimate
\[
s^2+s^2|\log s^2|
\le C_\varepsilon
\bigl(s^{2-\varepsilon}+s^{2+\varepsilon}\bigr),
\qquad s\ge0,
\]
combined with \eqref{eq03}, shows that the nonlinear
term in $E(u)$ is integrable for every
$u\in H^1(\R^N)$. Thus $E$ is well defined.
We next prove differentiability. Set
\[
G(z):=\frac12\bigl(|z|^2\log|z|^2-|z|^2\bigr),
\qquad
g(z):=z\log|z|^2,
\qquad z\in\C,
\]
with $G(0)=g(0)=0$. Identifying $\C$ with $\R^2$,
we have $G\in C^1(\R^2;\R)$ and $G^{\prime}(z)[w]=\RE\bigl(g(z)\overline w\bigr)$.
Let $\chi:[0,\infty)\to[0,1]$ be continuous, with
$\chi=1$ on $[0,1]$ and $\chi=0$ on $[2,\infty)$.
Decompose $g_-(z):=\chi(|z|)g(z)$ and $g_+(z):=(1-\chi(|z|))g(z)$.
These continuous functions satisfy
\[
|g_-(z)|\le C_\varepsilon|z|^{1-\varepsilon},
\qquad
|g_+(z)|\le C_\varepsilon|z|^{1+\varepsilon}.
\]
Consequently,
\begin{equation}
\label{Neo}
\begin{aligned}
\left\|
|x|^{-\frac{1-\varepsilon}{2-\varepsilon}}g_-(u)
\right\|_{L^{\frac{2-\varepsilon}{1-\varepsilon}}(\R^N)}
&\le C_\varepsilon
\left\|
|x|^{-\frac{1}{2-\varepsilon}}u
\right\|_{L^{2-\varepsilon}(\R^N)}^{1-\varepsilon},
\\
\left\|
|x|^{-\frac{1+\varepsilon}{2+\varepsilon}}g_+(u)
\right\|_{L^{\frac{2+\varepsilon}{1+\varepsilon}}(\R^N)}
&\le C_\varepsilon
\left\|
|x|^{-\frac{1}{2+\varepsilon}}u
\right\|_{L^{2+\varepsilon}(\R^N)}^{1+\varepsilon}.
\end{aligned}
\end{equation}
For $u\in H^1(\R^N)$, define the real linear
functional
\[
T(u)[h]
:=
\RE\int_{\R^N}\frac{g(u)\overline h}{|x|}\,dx.
\]
H\"older's inequality, \eqref{eq03} and \eqref{Neo} imply
\[
|T(u)[h]|
\le
C_{N,\varepsilon}
\left(
\|u\|_{H^1}^{1-\varepsilon}
+\|u\|_{H^1}^{1+\varepsilon}
\right)\|h\|_{H^1}.
\]
Moreover, combining \eqref{eq03} and \eqref{Neo}, and using the dominated convergence theorem, we obtain $T(u_n)\to T(u)$ in $H^{-1}$ whenever $u_n\to u$ in $H^1$.
On the other hand, define
\[
J(u):=\int_{\R^N}\frac{G(u)}{|x|}\,dx.
\]
The fundamental theorem of calculus gives
\[
J(u+h)-J(u)-T(u)[h]
=
\int_0^1\bigl(T(u+th)-T(u)\bigr)[h]\,dt.
\]
Thus, by  the continuity of $T$ at $u$ we obtain
\[
\begin{aligned}
|J(u+h)-J(u)-T(u)[h]|
&\le
\|h\|_{H^1}
\sup_{0\le t\le1}
\|T(u+th)-T(u)\|_{H^{-1}}
=o(\|h\|_{H^1})
\end{aligned}
\]
as $h\to0$ in $H^1$. Therefore $J$ is
Fr\'echet differentiable, with $J'(u)=T(u)$.
The continuity of $T$ proves that $J$ is of class
$C^1$. Finally, since $E(u)=\tfrac12\left\|\nabla u\right\|_{L^2(\R^N)}^2-J(u)$, we see that $E$ is of class $C^1$ with
$ E'(u)=-\Delta u-|x|^{-1}u\log|u|^2$.
\end{proof}

\begin{lemma}[Sharp  logarithmic Sobolev inequality]
\label{lem1}
Let $N\ge2$ and define
\[
\mathcal C_N
:=
\log\left[
\frac{\Gamma(N/2)}
{2\pi^{N/2}\Gamma(N-1)}
\left(\frac{4}{e(N-1)}\right)^{N-1}
\right].
\]
For every $u\in H^1(\R^N)\setminus\{0\}$, set $m_1(u):=\int_{\mathbb{R}^N}|x|^{-1}|u|^2\,dx$.
Then
\begin{equation}
\int_{\R^N}\frac{|u|^2}{|x|}
\log\left(\frac{|u|^2}{m_1(u)}\right)\,dx
\le
m_1(u)\left[
\mathcal C_N+(N-1)
\log\left(\frac{\left\|\nabla u\right\|_{L^2(\R^N)}^2}{m_1(u)}\right)
\right].
\label{eq04}
\end{equation}
The constant $\mathcal C_N$ is optimal. Equality holds
in \eqref{eq04} if and only if $u(x)=z e^{-a|x|}$  for a.e.  $x\in\R^N $
for some $a>0$ and $z\in\C\setminus\{0\}$.
\end{lemma}

\begin{proof}
For $N\ge 3$, the result follows immediately from the weighted logarithmic Sobolev inequality in \cite[Theorem~1.1 and Corollary~1.2]{DZ22}.
We consider $N=2$. In polar coordinates, define $v(s,\theta):=u(s^2,\theta)$.
The change of variables $r=s^2$ gives $m_1(u)=2\left\|v\right\|_{L^2(\R^2)}^2$
and
\begin{equation}
\begin{aligned}
\left\|\nabla u\right\|_{L^2(\R^2)}^2
&=
\frac12\left\|\nabla v\right\|_{L^2(\R^2)}^2
+\frac32\int_0^{2\pi}\int_0^\infty
\frac{|\partial_\theta v|^2}{s^2}\,s\,ds\,d\theta.
\end{aligned}
\label{eq05}
\end{equation}
In particular, $v\in H^1(\R^2)$ and $\left\|\nabla v\right\|_{L^2(\R^2)}^2\le2\left\|\nabla u\right\|_{L^2(\R^2)}^2$.
We recall the sharp euclidean logarithmic Sobolev inequality in dimension two (see \cite[Theorem
8.14]{LL01})
\begin{equation}
\label{LoS}
\int_{\mathbb{R}^2}|v|^2
\log\left(\frac{|v|^2}{\left\|v\right\|_{L^2(\R^2)}^2}\right)\,dy
\le
\left\|v\right\|_{L^2(\R^2)}^2
\log\left(
\frac{\left\|\nabla v\right\|_{L^2(\R^2)}^2}{\pi e\,\left\|v\right\|_{L^2(\R^2)}^2}
\right).
\end{equation}
On the other hand, the same change of variables yields
\[
\begin{aligned}
&\int_{\R^2}\frac{|u|^2}{|x|}
\log\left(\frac{|u|^2}{m_1(u)}\right)\,dx
=
2\int_{\R^2}|v|^2
\log\left(\frac{|v|^2}{\left\|v\right\|_{L^2(\R^2)}^2}\right)\,dy
-m_1(u)\log2.
\end{aligned}
\]
Combining these identities, we obtain
\begin{equation*}
\begin{aligned}
\int_{\mathbb{R}^2}\frac{|u|^2}{|x|}
\log\left(\frac{|u|^2}{m_1(u)}\right)\,dx
&\le m_1(u)\log\left(
\frac{\left\|\nabla v\right\|_{L^2(\R^2)}^2}{\pi e\,m_1(u)}
\right)
\le m_1(u)\log\left(
\frac{2\left\|\nabla u\right\|_{L^2(\R^2)}^2}{\pi e\,m_1(u)}
\right).
\end{aligned}
\end{equation*}
Since $\mathcal C_2=\log\left(\frac{2}{\pi e}\right)$, we obtain 
\eqref{eq04} for $N=2$. If equality holds, then equality must hold both
in the Euclidean logarithmic Sobolev inequality and in $\left\|\nabla v\right\|_{L^2(\R^2)}^2\le2\left\|\nabla u\right\|_{L^2(\R^2)}^2$.
By \eqref{eq05}, the latter forces
$\partial_\theta v=0$, so $v$ is radial about the
origin. The equality cases of the Euclidean
inequality then give $v(y)=z e^{-a|y|^2}$, where $a>0$ and $z\in\mathbb{C}\setminus\{0\}$.
Consequently, $u(x)=z e^{-a|x|}.$
Conversely, these profiles attain equality in all
the preceding inequalities. This proves optimality
and completes the classification of the equality
cases.
\end{proof}

\begin{lemma}\label{CCD}
Let $N\ge2$ and $0<\varepsilon<2/N$.
Suppose that $v_k\rightharpoonup v$ in $H^1(\R^N)$.
Then the following assertions hold.
\begin{enumerate}[label=\textup{(\roman*)}]
\item
For every $q\in\{2-\varepsilon,2+\varepsilon\}$,
\[
\left\||x|^{-1/q}(v_k-v)\right\|_{L^q(\R^N)}
\to 0.
\]
\item
We have
\[
\left\||x|^{-1/2}(v_k-v)\right\|_{L^2(\R^N)}
\to 0.
\]
\end{enumerate}
\end{lemma}

\begin{proof}
Fix $q\in\{2-\varepsilon,2+\varepsilon\}$.
The assumption $0<\varepsilon<2/N$ implies that
$qN/(N-1)>2$ and, when $N\ge3$,
$qN/(N-1)<2N/(N-2)$.
We may therefore choose $1<r<N$ sufficiently close
to $N$ so that, with $r'=r/(r-1)$, the exponent
$qr'$ belongs to the Sobolev embedding range.
For every $w\in H^1(\R^N)$, H\"older's inequality
and the Sobolev embedding then give
\begin{equation}
\int_{|x|<\delta}\frac{|w|^q}{|x|}\,dx
\le \left(\int_{|x|<\delta}|x|^{-r}\,dx\right)^{1/r}
\left(\int_{|x|<\delta}|w|^{qr'}\,dx\right)^{1/r'}
\le C\delta^{N/r-1}\|w\|_{H^1(\R^N)}^q.
\label{eqcd12}
\end{equation}
for $0<\delta\le1$. At infinity, the Sobolev embedding yields
\begin{equation}
\int_{|x|>R}\frac{|w|^{2+\varepsilon}}{|x|}\,dx
\le R^{-1}\left\|w\right\|_{L^{2+\varepsilon}(\R^N)}^{2+\varepsilon}
\le CR^{-1}\|w\|_{H^1(\R^N)}^{2+\varepsilon},
\qquad R\ge1.
\label{eqre13}
\end{equation}
Moreover, H\"older's inequality with
exponents $2/(2-\varepsilon)$ and $2/\varepsilon$ gives
\begin{equation}
\int_{|x|>R}\frac{|w|^{2-\varepsilon}}{|x|}\,dx
\le \left\|w\right\|_{L^2(\R^N)}^{2-\varepsilon}
\left(\int_{|x|>R}|x|^{-2/\varepsilon}\,dx\right)^{\varepsilon/2}
\le CR^{-1+N\varepsilon/2}\left\|w\right\|_{L^2(\R^N)}^{2-\varepsilon}.
\label{eqew14}
\end{equation}
The assumption $\varepsilon<2/N$ ensures both
the integrability of $|x|^{-2/\varepsilon}$ at
infinity and the decay of the last bound as
$R\to\infty$.

Now suppose that $v_k\rightharpoonup v$ in $H^1(\R^N)$,
and set $w_k=v_k-v$.
Then $M:=\sup_k\|w_k\|_{H^1(\R^N)}<\infty$.
For fixed $0<\delta<1<R$, let
$A_{\delta,R}:=\{x\in\R^N:\delta<|x|<R\}$.
By Rellich's theorem and uniqueness of the weak limit,
$w_k\to0$ in $L^q(A_{\delta,R})$.
When $q<2$, this convergence follows from convergence
in $L^2(A_{\delta,R})$, since the annulus has finite
measure. Hence
\[
\int_{A_{\delta,R}}\frac{|w_k|^q}{|x|}\,dx
\le \delta^{-1}\left\|w_k\right\|_{L^q(A_{\delta,R})}^q
\longrightarrow0.
\]
Combining this convergence with the estimates \eqref{eqcd12}--\eqref{eqew14}, we obtain
\[
\limsup_{k\to\infty}
\int_{\R^N}\frac{|w_k|^q}{|x|}\,dx
\le CM^q\bigl(\delta^{N/r-1}+R^{-\eta_q}\bigr),
\]
where $\eta_{2+\varepsilon}=1$ and
$\eta_{2-\varepsilon}=1-N\varepsilon/2>0$.
Letting $\delta\to0$ and $R\to\infty$ proves that
$\left\||x|^{-1/q}(v_k-v)\right\|_{L^q(\R^N)}\to0$.

Assertion~\textup{(ii)} follows from \textup{(i)}
and H\"older's inequality, since
\[
\begin{aligned}
\left\||x|^{-1/2}(v_k-v)\right\|_{L^2(\R^N)}^2
&\le
\left\||x|^{-\frac1{2-\varepsilon}}(v_k-v)\right\|_{L^{2-\varepsilon}(\R^N)}^{\frac{2-\varepsilon}{2}}
\left\||x|^{-\frac1{2+\varepsilon}}(v_k-v)\right\|_{L^{2+\varepsilon}(\R^N)}^{\frac{2+\varepsilon}{2}}
\to0.
\end{aligned}
\]
\end{proof}


\section{The Cauchy problem}
\label{caup}

In this section, we prove the existence and uniqueness
results stated in Propositions~\ref{pr01} and \ref{pr04}.
For existence, we regularize both the singular coefficient
and the logarithmic nonlinearity. Conservation of mass
and energy for the regularized problems provides uniform
bounds, while compactness in weighted Lebesgue spaces
allows us to pass to the limit. We then establish
uniqueness under the lower bound \eqref{un01}, using
a pointwise estimate for the difference of the nonlinear
terms and Strichartz estimates.

\begin{proof}[Proof of Proposition~\ref{pr01}]
For $n\ge2$, define (cf. \cite[Section~2]{CS21})
\[
b_n(x)=\frac{1}{\sqrt{|x|^2+n^{-2}}},
\qquad
\ell_n(s)=
\log\left(\frac{s+n^{-1}}{1+s/n}\right),
\qquad s\ge0,
\]
and $H_n(s)=\int_0^s\ell_n(\tau)\,d\tau$.
We consider the regularized Cauchy problem
\begin{equation}
i\partial_tu_n+\Delta u_n
+b_n(x)u_n\ell_n(|u_n|^2)=0,
\qquad u_n(0)=u_0.
\label{eq07}
\end{equation}
For each fixed $n$, we have
$b_n\in W^{1,\infty}(\R^N)$, and the map
$z\mapsto z\ell_n(|z|^2)$ is continuously differentiable,
vanishes at zero, and has bounded real derivative; see  \cite[Lemma~8]{CS21} for more details.
Applying \cite[Theorem~3.7.1]{C03}, we obtain a unique global solution
\[
u_n\in C(\mathbb{R};H^1(\mathbb{R}^N))
\cap C^1(\mathbb{R};H^{-1}(\mathbb{R}^N))
\]
such that both the mass and the energy are conserved, that is,
\begin{equation}
M(u_n(t))=M(u_0),
\qquad
E_n(u_n(t))=E_n(u_0),
\qquad t\in\mathbb{R},
\label{eq08}
\end{equation}
where
\[
E_n(v)=\frac12\left\|\nabla v\right\|_{L^2(\R^N)}^2
-\frac12\int_{\mathbb{R}^N}b_n(x)H_n(|v|^2)\,dx.
\]

We next establish estimates independent of $n$. 
For every $s>0$, $|\ell_n(s)|\le|\log s| $, $\ell_n(s)\to\log s $ 
and $\ell_n(s)$ has the same sign as $\log s$  (cf. \cite[Lemma~8]{CS21}).
Fix $0<\varepsilon<2/N$. It follows that
\begin{equation}
\begin{aligned}
 &|H_n(s)|\le C\bigl(s^{1-\varepsilon/2}+s^{1+\varepsilon/2}\bigr),
\qquad H_n(s)\le Cs^{1+\varepsilon/2},\\
&|z\ell_n(|z|^2)| \le C\bigl(|z|^{1-\varepsilon}+|z|^{1+\varepsilon}\bigr).
\end{aligned}
\label{eq09}
\end{equation}
uniformly in $n$. In addition, applying estimate \eqref{eq03} we obtain
\begin{equation}
\int_{\R^N}\frac{|v|^{2+\varepsilon}}{|x|}\,dx
\le
C\left\|\nabla v\right\|_{L^2(\R^N)}^{1+N\varepsilon/2}
\left\|v\right\|_{L^2(\R^N)}^{1-(N-2)\varepsilon/2},
\qquad v\in H^1(\R^N).
\label{eq06}
\end{equation}
Since $0<b_n(x)\le|x|^{-1}$, estimate \eqref{eq06}
and Young's inequality imply
\begin{equation}
E_n(v)\ge
\frac12\left\|\nabla v\right\|_{L^2(\R^N)}^2
-C\left\|\nabla v\right\|_{L^2(\R^N)}^{1+N\varepsilon/2}
\left\|v\right\|_{L^2(\R^N)}^{1-(N-2)\varepsilon/2}
\ge
\frac14\left\|\nabla v\right\|_{L^2(\R^N)}^2-C(N,\left\|v\right\|_{L^2(\R^N)}).
\label{Eax}
\end{equation}
where $C(N,\left\|v\right\|_{L^2(\R^N)})$ denotes a constant depending on $N$ and $\left\|v\right\|_{L^2(\R^N)}$.
Furthermore, $H_n(s)\to s\log s-s$. Thus, by \eqref{eq03}, \eqref{eq09}, and dominated convergence we see that $E_n(u_0)\to E(u_0)$, and $ \sup_{n\ge2}|E_n(u_0)|
\le C(N,\|u_0\|_{H^1})$.
Combining these facts with \eqref{eq08} and \eqref{Eax}, we obtain
\begin{equation}
\sup_{n\ge2}\sup_{t\in\R}
\|u_n(t)\|_{H^1(\R^N)}
\le C(N,\|u_0\|_{H^1}).
\label{eq10}
\end{equation}

The last estimate in \eqref{eq09}, H\"older's inequality, and \eqref{eq03} give
\begin{equation*}
\begin{aligned}
&\left|\int_{\R^N} b_n u_n\ell_n(|u_n|^2)\overline{h}\,dx\right|
\le C\sum_{q\in\{2-\varepsilon,\,2+\varepsilon\}}
\left(\int_{\R^N}\frac{|u_n|^q}{|x|}\,dx\right)^{(q-1)/q}
\left(\int_{\R^N}\frac{|h|^q}{|x|}\,dx\right)^{1/q}\\
&\qquad\le C\|h\|_{H^1(\R^N)}.
\end{aligned}
\end{equation*}

Hence \eqref{eq07} implies
$\sup_{n\ge2}\sup_{t\in\R}\|\partial_tu_n(t)\|_{H^{-1}}\le C$.
Interpolation between $H^{-1}(\R^N)$ and $H^1(\R^N)$ then yields
$\left\|u_n(t)-u_n(s)\right\|_{L^2(\R^N)}\le C|t-s|^{1/2}$ for all $s,t\in\R$.
By the Rellich and Arzel\`a--Ascoli theorems, followed by a
diagonal extraction, there exist a subsequence, still denoted by
$u_n$, and a function $u$ such that
$u_n\to u$ in $C([-T,T];L^2(B_R))$ for every $T,R>0$,
where $B_R=\{x\in\R^N:|x|<R\}$.
The uniform bounds also give
$u_n(t)\rightharpoonup u(t)$ in $H^1(\R^N)$ for every $t\in\R$,
and $u\in L^\infty(\R;H^1(\R^N))
\cap W_{\mathrm{loc}}^{1,\infty}(\R;H^{-1}(\R^N))$.
The boundedness in $H^1(\R^N)$ and continuity in
$H^{-1}(\R^N)$ imply that
$u\in C_w(\R;H^1(\R^N))$.
In particular, $u(0)=u_0$.

By Lemma~\ref{CCD}, the uniform $H^1$ bounds and
the local convergence imply that, for every $T>0$,
\begin{equation}
\sup_{|t|\le T}
\left\||x|^{-1/q}\bigl(u_n(t)-u(t)\bigr)\right\|_{L^q(\R^N)}
\to 0,
\qquad q\in\{2-\varepsilon,2+\varepsilon\}.
\label{eq12}
\end{equation}
Thus, we conclude that
$|x|^{-1/q}u_n\to |x|^{-1/q}u$
in $C([-T,T];L^q(\R^N))$ for $ q\in\{2-\varepsilon,2+\varepsilon\}$.

Since $ z\ell_n(|z|^2)\to z\log|z|^2$ locally uniformly in $z\in\C$
and $0<b_n(x)|x|\le1$ and
$b_n(x)|x|\to1$ almost everywhere,  and using \eqref{eq03},
we conclude that
\[
b_nu_n\ell_n(|u_n|^2)
\to |x|^{-1}u\log|u|^2
\quad\text{in }
L^1([-T,T];H^{-1}(\R^N)).
\]
Passing to the limit in \eqref{eq07} now proves that
$u$ satisfies \eqref{NLS} in the sense of distributions.
Since all terms belong to
$L^\infty_{\mathrm{loc}}(\R;H^{-1}(\R^N))$,
the equation holds in $H^{-1}(\R^N)$ for almost
every $t\in\R$. A standard argument shows that the mass is conserved.
To prove the energy inequality, fix $t\in\R$.
By \eqref{eq09}, \eqref{eq12}, we see that
\[
\int_{\R^N}b_n(x)H_n(|u_n(t)|^2)\,dx
\to
\int_{\R^N}
\frac{|u(t)|^2\log|u(t)|^2-|u(t)|^2}{|x|}\,dx.
\]
Weak lower semicontinuity  and
\eqref{eq08} therefore yield
\[
E(u(t))
\le\liminf_{n\to\infty}E_n(u_n(t))
=\lim_{n\to\infty}E_n(u_0)
=E(u_0).
\]
Since $t$ was arbitrary, this inequality holds for
every $t\in\R$. The asserted uniform $H^1$ bound
follows from \eqref{eq10} and weak lower
semicontinuity.

Finally, weak continuity in $H^1(\R^N)$ and Lemma \ref{CCD} imply
\[
u(t)\to u_0
\quad\text{in }
L^{2-\varepsilon}(\R^N,|x|^{-1}\,dx)
\cap
L^{2+\varepsilon}(\R^N,|x|^{-1}\,dx)
\]
as $t\to0$. Hence
\[
\int_{\R^N}
\frac{|u(t)|^2\log|u(t)|^2-|u(t)|^2}{|x|}\,dx
\to
\int_{\R^N}
\frac{|u_0|^2\log|u_0|^2-|u_0|^2}{|x|}\,dx.
\]
Together with $E(u(t))\le E(u_0)$, this gives $\limsup_{t\to0}\left\|\nabla u(t)\right\|_{L^2(\R^N)}^2
\le\left\|\nabla u_0\right\|_{L^2(\R^N)}^2$.

Since $u(t)\rightharpoonup u_0$ in $H^1(\R^N)$, we see that 
$\nabla u(t)\to\nabla u_0$ in $L^2(\R^N)$.
Combining this with continuity in $L^2(\R^N)$,
we obtain $\lim_{t\to0}\|u(t)-u_0\|_{H^1}=0$
This completes the proof.
\end{proof}

We next prove the uniqueness statement in
Proposition~\ref{pr04}. The lower bound \eqref{un01}
allows us to estimate the difference of the nonlinear
terms by a time-dependent potential, to which
Strichartz estimates can be applied.

\begin{proof}[Proof of Proposition~\ref{pr04}]
Set $w=u-v$ and $g(z)=z\log|z|^2$ with $g(0)=0$. We write
$\log_+s=\max\{\log s,0\}$ for $s>0$, with
$\log_+0=0$.

A direct computation yields the pointwise estimate (cf. \cite[Lemma 3.3]{F20})
\begin{equation}
|g(z)-g(\zeta)|
\le C|z-\zeta|
\bigl(1+|\log|\zeta||+\log_+|z|\bigr),
\qquad z\in\C,\quad \zeta\ne0.
\label{un03}
\end{equation}

Choose $\max\{1,N/2\}<p<N$.
We claim that
\begin{equation}
\frac{\log_+|u|+\log_+|v|}{|x|}
\in L^\infty(0,T;L^p).
\label{un04}
\end{equation}
Indeed, for every finite $s\ge1$
and every $f\in L^2$,
$\left\|\log_+|f|\right\|_{L^s(\R^N)}\le C_s\left\|f\right\|_{L^2(\R^N)}^{2/s}$,
because $(\log t)^s\le C_st^2$ for $t\ge1$.
Choose $p<\ell<N$ and $s=p\ell/(\ell-p)$.
Since $|x|^{-1}\in L^\ell(B_1)$,
H\"older's inequality gives
$\left\|\,|x|^{-1}\log_+|f|\,\right\|_{L^p(B_1)}
\le C\left\|\log_+|f|\right\|_{L^s(\R^N)}
\le C\left\|f\right\|_{L^2(\R^N)}^{2/s}$.
On $\R^N\setminus B_1$, we use $|x|^{-1}\le1$
to obtain
$\left\|\,|x|^{-1}\log_+|f|\,\right\|_{L^p(\R^N\setminus B_1)}
\le C\left\|f\right\|_{L^2(\R^N)}^{2/p}$.
Thus, conservation of mass implies that \eqref{un04} holds. By \eqref{un01},
$|\log|v(t,x)||\le
A_T(1+|x|)+\log_+|v(t,x)|$
almost everywhere. Consequently, \eqref{un03} and \eqref{un04}
implies
\begin{equation}
\left|\frac{g(u)-g(v)}{|x|}\right|
\le (C_T+V(t,x))|w|,
\qquad
V\in L^\infty(0,T;L^p),
\quad V\ge0,
\label{un05}
\end{equation}
where 
\[
V(t,x)=
C(1+A_T)\frac{\mathbf{1}_{B_1}(x)}{|x|}
+C\frac{\log_+|u(t,x)|+\log_+|v(t,x)|}{|x|}.
\]
Here $C_T$ is independent of time. Let $r=\frac{2p}{p-1}$ and $q=\frac{4p}{N}$.
Then $\frac{2}{q}+\frac{N}{r}=\frac{N}{2}$ and $q>2$,
so $(q,r)$ is an admissible pair. Moreover, $H^1\hookrightarrow L^r$.
Thus, as $w\in C([0,T];L^2)\cap L^\infty(0,T;H^1)$, we see that  on every interval $J=[t_0,t_1]\subset[0,T]$,
$w\in C(J;L^2)\cap L^q(J;L^r)$. For the application of the Strichartz
estimate, decompose $\frac{g(u)-g(v)}{|x|}=F_0+F_1$,
where
\[
F_0=
\frac{C_T}{C_T+V}\frac{g(u)-g(v)}{|x|},
\qquad
F_1=
\frac{V}{C_T+V}\frac{g(u)-g(v)}{|x|}.
\]
We may assume $C_T>0$.
By \eqref{un05}, we obtain $|F_0|\le C_T|w|$ and $ |F_1|\le V|w|$.
Since $1/r'=1/p+1/r$, H\"older's inequality
in space and time yields
\begin{equation}
\begin{aligned}
\left\|F_0\right\|_{L^1(J;L^2(\R^N))}
&\le C_T|J|\left\|w\right\|_{L^\infty(J;L^2(\R^N))},\\
\left\|F_1\right\|_{L^{q'}(J;L^{r'}(\R^N))}
&\le
\left\|V\right\|_{L^\infty(0,T;L^p(\R^N))}
|J|^{1-2/q}\left\|w\right\|_{L^q(J;L^r(\R^N))}.
\end{aligned}
\label{un06}
\end{equation}

Subtracting the equations for $u$ and $v$, we find $i\partial_tw+\Delta w=-(F_0+F_1)$
in the sense of distributions. We therefore have
\begin{equation}
w(t)=e^{i(t-t_0)\Delta}w(t_0)
+i\int_{t_0}^t
e^{i(t-s)\Delta}(F_0(s)+F_1(s))\,ds,
\qquad t\in J.
\label{un07}
\end{equation}

Applying
Strichartz estimates to \eqref{un07}, and using
\eqref{un06}, gives
\begin{equation}
\begin{aligned}
&\left\|w\right\|_{L^\infty(J;L^2(\R^N))}+\left\|w\right\|_{L^q(J;L^r(\R^N))}
\le C\left\|w(t_0)\right\|_{L^2(\R^N)}\\
&\quad+C\left(C_T|J|+\left\|V\right\|_{L^\infty(0,T;L^p(\R^N))}|J|^{1-2/q}\right)
\left(\left\|w\right\|_{L^\infty(J;L^2(\R^N))}+\left\|w\right\|_{L^q(J;L^r(\R^N))}\right).
\end{aligned}
\label{un08}
\end{equation}
Choose $\tau>0$ sufficiently small that $ C\left(
C_T\tau
+\left\|V\right\|_{L^\infty(0,T;L^p(\R^N))}\tau^{1-2/q}
\right)\le\frac12$.
For every $J$ of length at most $\tau$, \eqref{un08} yields $\left\|w\right\|_{L^\infty(J;L^2(\R^N))}
+\left\|w\right\|_{L^q(J;L^r(\R^N))}
\le 2C\left\|w(t_0)\right\|_{L^2(\R^N)}$.

Since $w(0)=0$, this proves that $w$ vanishes on
the first such interval. A finite partition of
$[0,T]$ into intervals of length at most $\tau$
therefore gives
$u(t)=v(t)$ in $L^2(\R^N)$ for all $t\in[0,T]$. Finally, since $u$ and $v  \in  C_w(\R;H^1(\R^N))$, \eqref{un02} holds in $H^1(\R^N)$ for every
$t\in[0,T]$, as claimed.
\end{proof}


\section{Standing waves and variational characterization}
\label{vari}

In this section, we prove Proposition~\ref{pr02} and
Theorem~\ref{th01}. We first construct explicit positive
solutions of the stationary problem \eqref{GS}.
We then use the sharp weighted logarithmic Sobolev
inequality to characterize all minimizers of \eqref{MPE}
for $\omega\ge0$ and show that the infimum is not
attained for $\omega<0$.

\begin{proof}[Proof of Proposition~\ref{pr02}]
We seek a positive solution of \eqref{GS} of the form
$\phi(x)=A e^{-a|x|}$, where $A>0$ and $a>0$.
For $x\ne0$, direct computation gives
\begin{equation}
\label{Lab}
-\Delta\phi+\omega\phi
-\frac{\phi\log|\phi|^2}{|x|}
=
\left(
\omega-a^2+2a
+\frac{(N-1)a-2\log A}{|x|}
\right)\phi.
\end{equation}
Since $\phi>0$, the stationary equation holds away
from the origin if and only if
$\omega=a^2-2a$ and $A=e^{\frac{N-1}{2}a}$.
The first identity is equivalent to $(a-1)^2=1+\omega$.
Thus real values of $a$ exist only for $\omega\ge -1$,
in which case $a=1+\sqrt{1+\omega}$ or
$a=1-\sqrt{1+\omega}$.
If $\omega=-1$, both roots equal $1$, giving
$\phi(x)=e^{\frac{N-1}{2}}e^{-|x|}$.
If $-1<\omega<0$, both roots are positive and distinct,
which yields the two profiles in assertion~\textup{(ii)}.
If $\omega\ge0$, only $1+\sqrt{1+\omega}$ is positive.
Finally, if $\omega<-1$, the quadratic equation
has no real root, so no positive solution of the
prescribed exponential form exists.
It remains to verify that the admissible profiles
solve \eqref{GS} on all of $\R^N$.
Note that for every $a>0$, the function $\phi(x)=e^{\frac{N-1}{2}a}e^{-a|x|}$
belongs to $H^1(\R^N)$. 
Since $\phi\in H^1(\R^N)$, Lemma~\ref{le01} implies that
\[
L:=-\Delta\phi+\omega\phi
-\frac{\phi\log|\phi|^2}{|x|}
\in H^{-1}(\R^N).
\]
Then, by \eqref{Lab}, $L$ vanishes on
$\mathbb{R}^N\setminus\{0\}$ with $N\geq 2$, so
$\operatorname{supp}L\subset\{0\}$, which implies
that $L=0$. Hence each admissible profile satisfies \eqref{GS}
in $H^{-1}(\R^N)$, completing the proof.
\end{proof}

We now turn to the minimization problem \eqref{MPE}.
The sharp inequality in Lemma~\ref{lem1} provides
a lower bound for the action on the Nehari manifold
and identifies all minimizers when $\omega\ge0$.

\begin{proof}[Proof of Theorem~\ref{th01}]
Recall that $m_1(u)=\int_{\mathbb{R}^N}|x|^{-1}|u|^2\,dx.$
By the definitions of $S_\omega$ and $I_\omega$, we see that
\begin{equation}
S_\omega(u)
=\frac12 I_\omega(u)+\frac12 m_1(u).
\label{ne01}
\end{equation}
In particular, $S_\omega(u)=\frac12m_1(u)>0$ for every nonzero $u\in H^1(\R^N)$ satisfying
$I_\omega(u)=0$. Suppose first that $\omega\ge0$.
Integration by parts and the Cauchy--Schwarz
inequality give
\[
(N-1)m_1(u)
=-2\operatorname{Re}\int_{\R^N}
\overline{u}\,\frac{x}{|x|}\cdot\nabla u\,dx
\le2\left\|u\right\|_{L^2(\R^N)}\left\|\nabla u\right\|_{L^2(\R^N)}.
\]
Consequently, for every nonzero $u\in H^1(\R^N)$,
\begin{equation}
\left\|u\right\|_{L^2(\R^N)}^2
\ge
\frac{(N-1)^2m_1(u)^2}{4\left\|\nabla u\right\|_{L^2(\R^N)}^2}.
\label{ne02}
\end{equation}

Let $u\ne0$ satisfy $I_\omega(u)=0$. We claim that $m_1(u)\ge m_1(\phi)$, where $\phi(x)=
e^{\frac{N-1}{2}(1+\sqrt{1+\omega})}
e^{-(1+\sqrt{1+\omega})|x|}$ is the profile given in Proposition \ref{pr02}(iii). Indeed,
by Lemma~\ref{lem1},
\[
\begin{aligned}
\left\|\nabla u\right\|_{L^2(\R^N)}^2+\omega\left\|u\right\|_{L^2(\R^N)}^2
&\le
m_1(u)\left[
\log m_1(u)+\mathcal C_N
+(N-1)\log\left(
\frac{\left\|\nabla u\right\|_{L^2(\R^N)}^2}{m_1(u)}
\right)
\right].
\end{aligned}
\]
Since $\omega\ge0$, we may use \eqref{ne02}
to obtain
\begin{equation}
\log m_1(u)\ge
\frac{\left\|\nabla u\right\|_{L^2(\R^N)}^2}{m_1(u)}
+\frac{\omega(N-1)^2m_1(u)}{4\left\|\nabla u\right\|_{L^2(\R^N)}^2}
-(N-1)\log\left(\frac{\left\|\nabla u\right\|_{L^2(\R^N)}^2}{m_1(u)}\right)
-\mathcal C_N.
\label{ne03}
\end{equation}
Notice also that $I_\omega(\phi)=0$,
and
\begin{equation}
\label{tios}
\frac{\left\|\nabla\phi\right\|_{L^2(\R^N)}^2}{m_1(\phi)}
=\frac{N-1}{2}(1+\sqrt{1+\omega}),
\qquad
\frac{\left\|\phi\right\|_{L^2(\R^N)}^2}{m_1(\phi)}
=\frac{N-1}{2(1+\sqrt{1+\omega})}.
\end{equation}
Now, consider the function (cf. \eqref{ne03})
\[
f(t)=
t+\frac{\omega(N-1)^2}{4t}
-(N-1)\log t-\mathcal C_N,
\qquad t>0.
\]
A direct calculation shows that it has a unique global minimum 
at $t_*=\frac{N-1}{2}(1+\sqrt{1+\omega})$.
Thus $\phi$ attains equality in \eqref{ne02}
and realizes the scalar minimum above.
It also attains equality in \eqref{eq04},
by the equality characterization in
Lemma~\ref{lem1}.
Hence
\[
\log m_1(u)
\ge f\left(\frac{\left\|\nabla u\right\|_{L^2(\R^N)}^2}{m_1(u)}\right)
\ge f(t_*)
=\log m_1(\phi),
\]
where $t_*=\frac{N-1}{2}(1+\sqrt{1+\omega})$.
Therefore $m_1(u)\ge m_1(\phi)$, as claimed.
Moreover, by \eqref{ne01}, we see that $d(\omega)=S_\omega(\phi)=\frac12m_1(\phi)$.
In particular,
we obtain
\[
d(\omega)=
\frac{|\mathbb S^{N-1}|\Gamma(N-1)}
{2[2(1+\sqrt{1+\omega})]^{N-1}}
e^{(N-1)(1+\sqrt{1+\omega})},
\]
where $\Gamma$ denotes the gamma function.
To characterize all minimizers, let
$u\in\mathcal N_\omega$.
Then $m_1(u)=m_1(\phi)$, so equality must hold
both in the logarithmic Sobolev inequality and
in the scalar minimization used above.
Lemma~\ref{lem1} therefore gives
$u(x)=z e^{-a|x|}$, where $a>0$ and $z\in\mathbb{C}\setminus\{0\}$.
For such a function,
\[
\frac{\left\|\nabla u\right\|_{L^2(\R^N)}^2}{m_1(u)}
=\frac{N-1}{2}a.
\]
Uniqueness of the scalar minimum $t_*$ and \eqref{tios} implies that $ a=1+\sqrt{1+\omega}$.
On the other hand, direct integration gives
\[
\frac{\left\|u\right\|_{L^2(\R^N)}^2}{m_1(u)}
=\frac{N-1}{2a},
\qquad
\frac{1}{m_1(u)}
\int_{\R^N}\frac{|u|^2\log|u|^2}{|x|}\,dx
=\log|z|^2-(N-1).
\]
The constraint $I_\omega(u)=0$ consequently yields
$ \log|z|^2
=N-1+\tfrac{N-1}{2}
\left(a+\tfrac{\omega}{a}\right).$
Since $\omega=a^2-2a$, this reduces to $ |z|^2=e^{(N-1)a}$.
Therefore $u=e^{i\theta}\phi$ for some $\theta\in\R$.
Conversely, every phase rotation of $\phi$ is
a minimizer, since both $S_\omega$ and $I_\omega$
are invariant under phase rotations.
This proves assertion~\textup{(i)}, including
the case $\omega=0$.

Suppose now that $\omega<0$.
For $a>0$, consider
\[
v_a(x)=
e^{\frac{N-1}{2}
+\frac{N-1}{4}\left(a+\frac{\omega}{a}\right)}
e^{-a|x|}.
\]
Each $v_a \in H^1(\R^N)$. Moreover, a direct calculation shows that 
$I_\omega(v_a)=0$ for every $a>0$.
Using \eqref{ne01}, we obtain
\[
S_\omega(v_a)
=
\frac{|\mathbb S^{N-1}|\Gamma(N-1)}
{2(2a)^{N-1}}
e^{N-1+\frac{N-1}{2}
\left(a+\frac{\omega}{a}\right)}
\to0
\qquad\text{as }a\to0^+.
\]
It follows that $d(\omega)=0$.
Finally, \eqref{ne01} shows that the action is
strictly positive at every nonzero function
satisfying the Nehari constraint.
Hence the infimum is not attained, proving
assertion~\textup{(ii)}.
\end{proof}


\section{Stability of the ground states}
\label{stab}

In this section, we prove the orbital stability result
stated in Theorem~\ref{th02}. The argument relies on
the variational characterization in Theorem~\ref{th01},
the compactness of minimizing sequences, and the strict
convexity of $d(\omega)$ for $\omega>0$. Conservation
of mass and the energy inequality then ensure that
solutions starting near the phase orbit remain close
to it for all time.

\begin{proof}[Proof of Theorem~\ref{th02}]

Fix $\omega>0$ and write $\phi(x)=e^{\frac{N-1}{2}(1+\sqrt{1+\omega})}e^{-(1+\sqrt{1+\omega})|x|}$.
By Theorem~\ref{th01}, $\mathcal N_\omega=\{e^{i\theta}\phi:\theta\in\R\}$ and $S_\omega(\phi)=d(\omega)$.
For every $\nu>0$, identity \eqref{ne01} gives $S_\nu(v)=\frac12I_\nu(v)+\frac12m_1(v)$.

We claim that if $ \frac12m_1(v)=d(\nu)$, then
\begin{equation}
I_\nu(v)\ge0
\quad\text{and}\quad
S_\nu(v)\ge d(\nu).
\label{st01}
\end{equation}
Indeed, for every $\lambda>0$, $I_\nu(\lambda v)=\lambda^2\bigl(I_\nu(v)-m_1(v)\log\lambda^2\bigr)$.
If $I_\nu(v)<0$, choosing $\lambda=e^{I_\nu(v)/(2m_1(v))}<1$ gives $I_\nu(\lambda v)=0$, and hence $d(\nu)\le S_\nu(\lambda v)=\frac{\lambda^2}{2}m_1(v)=\lambda^2d(\nu)<d(\nu)$, a contradiction. This proves the claim.

Next, differentiating the explicit formula for $d$ in
Theorem~\ref{th01}, we obtain, for $\nu>0$,
\[
d'(\nu)=
\frac{|\mathbb S^{N-1}|\Gamma(N)}
{2[2(1+\sqrt{1+\nu})]^N}
e^{(N-1)(1+\sqrt{1+\nu})}>0
\]
and
\[
d''(\nu)
=
d'(\nu)
\frac{(N-1)(1+\sqrt{1+\nu})-N}
{2(1+\sqrt{1+\nu})\sqrt{1+\nu}}
>0.
\]
In particular, $ d'(\omega)=\frac12\left\|\phi\right\|_{L^2(\R^N)}^2$.
Choose $0<\eta<\omega$ sufficiently small that
$d''(\nu)\ge\frac12d''(\omega)$ for every
$\nu\in[\omega-\eta,\omega+\eta]$.
Taylor's formula then gives
\begin{equation}
d(\nu)-d(\omega)-d'(\omega)(\nu-\omega)
\ge\frac14d''(\omega)|\nu-\omega|^2
\label{st02}
\end{equation}
on this interval. Since $d$ is strictly increasing, whenever
$d(\omega-\eta)\le\frac12m_1(v)\le d(\omega+\eta)$,
there exists a unique $\nu(v)\in[\omega-\eta,\omega+\eta]$
such that $d(\nu(v))=\frac12m_1(v)$.
Therefore, combining \eqref{st01} and \eqref{st02}, we obtain
\begin{equation}
\begin{aligned}
&E(v)-E(\phi)
+\frac{\nu(v)}2
\bigl(\left\|v\right\|_{L^2(\R^N)}^2-\left\|\phi\right\|_{L^2(\R^N)}^2\bigr)
= S_{\nu(v)}(v)-d(\omega)
-\frac{\nu(v)-\omega}{2}\left\|\phi\right\|_{L^2(\R^N)}^2\\
&\ge d(\nu(v))-d(\omega)
-d'(\omega)(\nu(v)-\omega)
\ge \frac14d''(\omega)|\nu(v)-\omega|^2.
\end{aligned}
\label{st03}
\end{equation}

We next prove that every bounded sequence
$(v_n)$ in $H^1(\R^N)$ satisfying
\begin{equation}
S_\omega(v_n)\to d(\omega),
\qquad
I_\omega(v_n)\to0
\label{st04}
\end{equation}
also satisfies
\begin{equation}
\inf_{\theta\in\R}
\|v_n-e^{i\theta}\phi\|_{H^1}\to0.
\label{st05}
\end{equation}
After extracting a subsequence, we may assume that
$v_n\rightharpoonup v$ in $H^1(\R^N)$.
Lemma~\ref{CCD}, interpolation, and the logarithmic
growth estimates in Lemma~\ref{le01}, yield
\begin{equation}
\begin{aligned}
m_1(v_n)&\to m_1(v),\\
\int_{\mathbb{R}^N}\frac{|v_n|^2\log|v_n|^2}{|x|}\,dx
&\to
\int_{\mathbb{R}^N}\frac{|v|^2\log|v|^2}{|x|}\,dx.
\end{aligned}
\label{lbnr}
\end{equation}
By \eqref{ne01} and \eqref{st04}, $m_1(v)=2d(\omega)>0$, so $v\ne0$. Since $\omega>0$, weak lower semicontinuity and \eqref{lbnr} gives $I_\omega(v)\le\liminf_{n\to\infty}I_\omega(v_n)=0$. On the other hand, \eqref{st01} implies $I_\omega(v)\ge0$. Thus $I_\omega(v)=0$ and $S_\omega(v)=d(\omega)$. Theorem~\ref{th01} therefore yields $v=e^{i\theta}\phi$ for some $\theta\in\R$. Moreover, by \eqref{lbnr} we infer that $\left\|\nabla v_n\right\|_{L^2(\R^N)}^2+\omega\left\|v_n\right\|_{L^2(\R^N)}^2
\to
\left\|\nabla v\right\|_{L^2(\R^N)}^2+\omega\left\|v\right\|_{L^2(\R^N)}^2$,  together with weak convergence,
it implies $v_n\to v$ strongly in $H^1(\R^N)$.
The same argument applies to every subsequence,
proving \eqref{st05} for the entire sequence.

We now argue by contradiction.
If the standing wave were not orbitally stable,
there would exist $\varepsilon_0>0$, global weak
solutions $u_n$ satisfying the properties in
Proposition~\ref{pr01}, and times $t_n\in\R$
such that $\|u_n(0)-\phi\|_{H^1}\to0$
whereas
\begin{equation}
\inf_{\theta\in\R}
\|u_n(t_n)-e^{i\theta}\phi\|_{H^1}
\ge\varepsilon_0.
\label{st06}
\end{equation}
By continuity of $E$ and $M$ on $H^1(\mathbb{R}^N)$, $E(u_n(0))\to E(\phi)$ and $M(u_n(0))\to M(\phi)$.

The coercive estimate (cf. \eqref{Eax}) used in the proof of
Proposition~\ref{pr01} gives $E(v)\ge\frac14\left\|\nabla v\right\|_{L^2(\R^N)}^2-C$. The energy inequality, $E(u(t))\le E(u_0)$ for all $t\in\R$,  imply
\begin{equation}
\sup_n\sup_{t\in\R}\|u_n(t)\|_{H^1}<\infty.
\label{st07}
\end{equation}
Furthermore, weak continuity in $H^1(\R^N)$
and compactness of the weighted embeddings (cf. Lemma \ref{CCD}) imply
that $t\mapsto m_1(u_n(t))$
is continuous for every $n$.

We claim that, for all sufficiently large $n$,
\begin{equation}
d(\omega-\eta)<\frac12m_1(u_n(t))<d(\omega+\eta)
\qquad\text{for every }t\in\R.
\label{st08}
\end{equation}
Since $d$ is strictly increasing and
$m_1(u_n(0))/2\to d(\omega)$, these inequalities
hold at $t=0$ for all sufficiently large $n$.
Whenever $m_1(u_n(t))/2$ belongs to the corresponding
closed interval, $\nu(u_n(t))\in[\omega-\eta,\omega+\eta]$.
Thus \eqref{st03}, mass conservation, and the energy
inequality give
\begin{equation}
\frac14d''(\omega)|\nu(u_n(t))-\omega|^2
\le |E(u_n(0))-E(\phi)|
+\frac{\omega+\eta}{2}|M(u_n(0))-M(\phi)|.
\label{st09}
\end{equation}
The right-hand side is independent of $t$ and tends
to zero. For all sufficiently large $n$, it is
therefore strictly smaller than
$d''(\omega)\eta^2/4$.
If $m_1(u_n(t))/2$ reached either endpoint, the
definition of $\nu$ would give
$\nu(u_n(t))=\omega\pm\eta$, contradicting
\eqref{st09}. Since $t\mapsto m_1(u_n(t))$ is
continuous and starts inside the interval, it
cannot leave it in either time direction.
This proves \eqref{st08}.

It now follows from \eqref{st09} that
$\sup_{t\in\mathbb{R}}|\nu(u_n(t))-\omega|\to0$.
Consequently, $m_1(u_n(t_n))=2d(\nu(u_n(t_n)))\to2d(\omega)$.
In addition, \eqref{st03}, \eqref{st08} and mass conservation
give
\[
E(u_n(t_n))-E(\phi)
\ge
-\frac{\omega+\eta}{2}
|M(u_n(0))-M(\phi)|.
\]

Combining this lower bound with $E(u_n(t_n))\le E(u_n(0))\to E(\phi)$, we obtain $E(u_n(t_n))\to E(\phi)$. Mass conservation and \eqref{ne01} then imply $S_\omega(u_n(t_n))\to d(\omega)$ and $I_\omega(u_n(t_n))\to0$. By \eqref{st07}, the sequence $(u_n(t_n))$ is bounded in $H^1(\mathbb{R}^N)$. The compactness property \eqref{st05} therefore gives $\inf_{\theta\in\mathbb{R}}\|u_n(t_n)-e^{i\theta}\phi\|_{H^1}\to0$, contradicting \eqref{st06}. This proves the orbital stability of $e^{i\omega t}\phi$.
\end{proof}




\section{Orbital instability of standing waves}
\label{inst}


The spectral properties and coercivity estimates
used below are established in Appendix~\ref{asp}.
We fix the functions $\xi,d,\eta$ and the constants
$\lambda,\kappa$ supplied by Lemma~\ref{sp02},
with the normalization $(\xi,d)_{L^2}=1$.

We define $g(z)=z\log|z|^2$ and $G(z)=\frac12\bigl(|z|^2\log|z|^2-|z|^2\bigr)$
with $g(0)=G(0)=0$.
For $\rho>0$ and $z\in\C$, let
\begin{equation}
\begin{aligned}
\mathcal R(\rho,z)
={}&G(\rho+z)-G(\rho)
-\rho\log(\rho^2)\operatorname{Re}z\\
&-\frac12(\log\rho^2+2)(\operatorname{Re}z)^2
-\frac12\log(\rho^2)(\operatorname{Im}z)^2
\end{aligned}
\label{ErF}
\end{equation}
and
\begin{equation}
\label{Nonl}
\mathcal N(\rho,z)
=
g(\rho+z)-g(\rho)
-\log(\rho^2)z-2\operatorname{Re}z.
\end{equation}

\begin{lemma}[Estimates for the logarithmic remainders]
\label{le03}
Fix $0<\alpha<2/N$.
There exist constants $C,C_\alpha>0$ such that,
for every $\rho>0$ and $z\in\C$,
\begin{equation}
\label{Rees}
|\mathcal R(\rho,z)|\le C_\alpha\rho^{-\alpha}|z|^{2+\alpha},
\qquad
|\mathcal R(\rho,z)|\le C|z|^2
\left(1+\log_+\frac{|z|}{\rho}\right).
\end{equation}
and
\begin{equation}
\label{Nes}
|\mathcal N(\rho,z)|
\le C_\alpha\rho^{-\alpha}|z|^{1+\alpha}.
\end{equation}
For $R\ge1$, let $\chi_R$ be a measurable function
satisfying $0\le\chi_R\le1$, $\chi_R=0$ on $\{|x|\le R\}$,
and $\chi_R=1$ on $\{|x|\ge2R\}$ and set $Y_R(h)=\int_{\R^N}\chi_R|h|^2\,dx$.

There exist $C>0$, independent of $R$ and of the
choice of $\chi_R$, and $C_R>0$ such that
\begin{equation}
\label{NLyu}
\int_{\mathbb{R}^N}
\frac{|\mathcal R(\phi,h)|}{|x|}\,dx
\le
C Y_R(h)+C_R\|h\|_{H^1}^{2+\alpha},
\qquad h\in H^1(\mathbb{R}^N).
\end{equation}
Moreover,
\begin{equation}
\label{eqn1}
\int_{\mathbb{R}^N}
\frac{|\mathcal N(\phi,h)|\,|\psi|}{|x|}\,dx
\le C\|h\|_{H^1}^{1+\alpha},
\qquad
\psi\in\{\phi,\xi,d\}.
\end{equation}
\end{lemma}
\begin{proof}
Fix $0<\alpha<2/N$. Since $N\ge2$, we have
$0<\alpha<1$. We first establish the estimates \eqref{Rees} and \eqref{Nes}.
The definitions of $G$ and $g$ give the scaling
identities
\begin{equation}
\label{scal}
\mathcal R(\rho,z)
=\rho^2\mathcal R\left(1,\tfrac z\rho\right),
\qquad
\mathcal N(\rho,z)
=\rho\,\mathcal N\left(1,\tfrac z\rho\right).
\end{equation}
It therefore suffices to
estimate them when $\rho=1$.
Since $G$ and $g$ are smooth in a neighborhood
of $1$, Taylor's formula gives $\mathcal R(1,z)=O(|z|^3)$ and
$\mathcal N(1,z)=O(|z|^2)$ as $z\to0$.
Both functions $z\mapsto\mathcal R(1,z)$ and
$z\mapsto\mathcal N(1,z)$ are continuous on $\mathbb{C}$, including
at $z=-1$, because $G$ and $g$ are continuous
at zero. Moreover, their explicit formulas imply
$|\mathcal R(1,z)|\le C|z|^2(1+\log|z|)$ and
$|\mathcal N(1,z)|\le C|z|(1+\log|z|)$ for $|z|\ge2$.
Combining the estimates near zero, continuity
on compact sets away from zero, and the bounds
at infinity, we obtain $|\mathcal R(1,z)|\le C_\alpha|z|^{2+\alpha}$,
$|\mathcal R(1,z)|\le C|z|^2(1+\log_+|z|)$, and
$|\mathcal N(1,z)|\le C_\alpha|z|^{1+\alpha}$.
Here we used $\alpha<1$ near zero and
$\log s\le C_\alpha s^\alpha$ for $s\ge1$.
The scaling identities \eqref{scal} now yield the estimates \eqref{Rees} and \eqref{Nes}.

We next estimate \eqref{NLyu}. On $\{|x|\le2R\}$, the explicit profile satisfies
\[
\phi(x)\ge
e^{\frac{N-1}{2}(1+\sqrt{1+\omega})
-2R(1+\sqrt{1+\omega})}>0.
\]
The first scalar estimate in \eqref{Rees} and \eqref{eq03}
therefore give
\begin{equation}
\label{eqr2}
\int_{|x|\le2R}
\frac{|\mathcal R(\phi,h)|}{|x|}\,dx
\le
C_R\int_{\mathbb R^N}\frac{|h|^{2+\alpha}}{|x|}\,dx
\le C_R\|h\|_{H^1}^{2+\alpha}.
\end{equation}
On $\{|x|>2R\}$, we use
\[
\log_+\frac{|h|}{\phi}
\le\log_+|h|+\log_+\frac1\phi
\le
\log_+|h|+(1+\sqrt{1+\omega})|x|.
\]
Since $R\ge1$, the second scalar estimate in \eqref{Rees} implies
\[
\frac{|\mathcal R(\phi,h)|}{|x|}
\le
C|h|^2
+C_\alpha\frac{|h|^{2+\alpha}}{|x|}
\qquad\text{on }\{|x|>2R\}.
\]
Using $\chi_R=1$ on this region and
\eqref{eq03}, we obtain
\begin{equation}
\label{Red}
\int_{|x|>2R}
\frac{|\mathcal R(\phi,h)|}{|x|}\,dx
\le
C Y_R(h)+C\|h\|_{H^1}^{2+\alpha}.
\end{equation}
where $C$ is independent of $R$ and of the
choice of $\chi_R$. Combining \eqref{eqr2} and \eqref{Red} we obtain \eqref{NLyu}.
Next, we prove \eqref{eqn1}. Let
$\psi\in\{\phi,\xi,d\}$.
By Lemma~\ref{le02},
$\phi^{-\alpha}\psi\in H^1(\R^N)$.
The estimate for $\mathcal N$  (cf. \eqref{Nes})  and
H\"older's inequality  yield
\[
\begin{aligned}
\int_{\R^N}
\frac{|\mathcal N(\phi,h)|\,|\psi|}{|x|}\,dx
&\le
C_\alpha\int_{\R^N}
\frac{|h|^{1+\alpha}|\phi^{-\alpha}\psi|}{|x|}\,dx\\
&\le
C_\alpha
\left(
\int_{\R^N}\frac{|h|^{2+\alpha}}{|x|}\,dx
\right)^{\frac{1+\alpha}{2+\alpha}}
\left(
\int_{\R^N}
\frac{|\phi^{-\alpha}\psi|^{2+\alpha}}{|x|}\,dx
\right)^{\frac1{2+\alpha}}\\
&\le C\|h\|_{H^1}^{1+\alpha},
\end{aligned}
\]
where the last inequality follows from
\eqref{eq03}, which implies \eqref{eqn1}.

\end{proof}

\begin{lemma}\label{le07}
For every $h=u+iv\in H^1_{\mathrm{rad}}(\R^N)$,
with $u,v$ real-valued, we have
\[
S_\omega(\phi+h)-S_\omega(\phi)
=
\frac12\mathcal F_+[u]
+\frac12\mathcal F_-[v]
-\int_{\R^N}\frac{\mathcal R(\phi,h)}{|x|}\,dx.
\]
\end{lemma}
\begin{proof}
Let $h=u+iv\in H^1_{\mathrm{rad}}(\R^N)$,
with $u,v$ real-valued.
Substituting the definition of $\mathcal R$
into the action gives
\begin{align}
S_\omega(\phi+h)-S_\omega(\phi)
={}&
\int_{\R^N}\nabla\phi\cdot\nabla u\,dx
+\omega\int_{\R^N}\phi u\,dx
-\int_{\R^N}\frac{\phi\log\phi^2}{|x|}u\,dx
\label{ac01}\\
&+\frac12\int_{\R^N}
\bigl(|\nabla u|^2+|\nabla v|^2
+\omega(u^2+v^2)\bigr)\,dx
\label{ac02}\\
&-\frac12\int_{\R^N}
\frac{(\log\phi^2+2)u^2+\log\phi^2v^2}{|x|}\,dx
\label{ac03}\\
&-\int_{\R^N}\frac{\mathcal R(\phi,h)}{|x|}\,dx.
\label{Rfa}
\end{align}
The remainder\eqref{Rfa} integrable by \eqref{NLyu}.
The terms on the right-hand side of \eqref{ac01}
sum to zero by the stationary equation
\eqref{GS}. Moreover, since
$\log\phi^2=(N-1)a-2a|x|$ and $\omega+2a=a^2$,
the sum of the terms in \eqref{ac02} and
\eqref{ac03} equals
$\frac12\mathcal F_+[u]+\frac12\mathcal F_-[v]$.
This proves the asserted identity.
\end{proof}

For $h=u+iv\in H^1_{\mathrm{rad}}(\R^N)$,
with $u,v$ real-valued, define
\[
z_\pm(h)=\lambda(u,d)_{L^2}\pm(v,\xi)_{L^2}.
\]

\begin{lemma}[Localized energy coercivity]
\label{le04}
There exists $C_*>0$ with the following property. For every
$R\ge1$, there exists $\delta_R>0$ such that, if
$h=u+iv\in H^1_{\mathrm{rad}}(\mathbb R^N)$ satisfies
$\|h\|_{H^1}\le\delta_R$ and $(v,\phi)_{L^2}=0$, and
$M(\phi+h)=M(\phi)$ and $E(\phi+h)\le E(\phi)$, then
\begin{equation}
\label{estih}
\|h\|_{H^1}^2
\le C_*\bigl((u,d)_{L^2}^2+Y_R(h)\bigr).
\end{equation}
After increasing $C_*$ if necessary, we also have
\begin{equation}\label{vity}
\|h\|_{H^1}^2\le C_*\bigl(|z_+(h)|^2+|z_-(h)|^2+Y_R(h)\bigr).
\end{equation}
Here $\chi_R$ and $Y_R$ are as in Lemma~\ref{le03}.
The constant $C_*$ is independent of $R$ and of
the choice of $\chi_R$.
\end{lemma}
\begin{proof}
Fix $0<\alpha<2/N$, and let $h=u+iv$ satisfy
the hypotheses of the lemma.
The mass constraint gives $2(u,\phi)_{L^2}+\left\|h\right\|_{L^2(\R^N)}^2=0$.
Decompose
\[
u=t\xi+u_\perp+\mu\phi,
\qquad
t=(u,d)_{L^2},
\qquad
\mu=-\frac{\left\|h\right\|_{L^2(\R^N)}^2}{2M(\phi)}.
\]
By Lemma~\ref{sp02}, $\xi$ and $d$ are orthogonal
to $\phi$ and $(\xi,d)_{L^2}=1$.
Consequently, $(u_\perp,\phi)_{L^2}=(u_\perp,d)_{L^2}=0$.
Since $(v,\phi)_{L^2}=0$, Lemmas~\ref{sp01}
and~\ref{sp03} yield a constant $k>0$ such that
\begin{equation}
\label{eqfg}
\mathcal F_+[u_\perp]\ge k\|u_\perp\|_{H^1}^2,
\qquad
\mathcal F_-[v]\ge k\|v\|_{H^1}^2.
\end{equation}
Also, the spectral identities \eqref{df} imply
\begin{align}
&\mathcal F_+[\xi]=-\lambda^2,
\qquad
\mathcal F_+(\xi,u_\perp)=0,
\label{igh}\\
&\mathcal F_+[t\xi+u_\perp]
=-\lambda^2t^2+\mathcal F_+[u_\perp].
\label{uji}
\end{align}
Moreover, $|\mu|\le C\|h\|_{H^1}^2$ and
$t\xi+u_\perp=u-\mu\phi$.
By continuity of the bilinear form $\mathcal F_+$
on $H^1_{\mathrm{rad}}\times H^1_{\mathrm{rad}}$,
$\left|
\mathcal F_+[u]-\mathcal F_+[t\xi+u_\perp]
\right|
\le C\|h\|_{H^1}^3$
whenever $\|h\|_{H^1}\le1$. By \eqref{uji} and \eqref{eqfg},
it follows that
\begin{equation}\label{co02}
\mathcal F_+[u]+\mathcal F_-[v]
\ge
-\lambda^2t^2
+k\bigl(\|u_\perp\|_{H^1}^2+\|v\|_{H^1}^2\bigr)
-C\|h\|_{H^1}^3.
\end{equation}

The decomposition of $u$ also gives
\[
\|h\|_{H^1}^2
\le
C_0\bigl(t^2+\|u_\perp\|_{H^1}^2+\|v\|_{H^1}^2\bigr)
+C_0\|h\|_{H^1}^4.
\]
Taking $\|h\|_{H^1}$ below a fixed threshold,
independent of $R$, we may absorb the last term.
Combining the resulting estimate with
\eqref{co02}, we obtain constants $k_0,C>0$,
independent of $R$, such that
\begin{equation}\label{co03}
\frac12\mathcal F_+[u]+\frac12\mathcal F_-[v]
\ge
k_0\|h\|_{H^1}^2-Ct^2-C\|h\|_{H^1}^3.
\end{equation}

Since $M(\phi+h)=M(\phi)$, Lemma~\ref{le07}
gives
\[
E(\phi+h)-E(\phi)
=
\frac12\mathcal F_+[u]+\frac12\mathcal F_-[v]
-\int_{\R^N}\frac{\mathcal R(\phi,h)}{|x|}\,dx.
\]
Using \eqref{co03} and the 
estimate \eqref{NLyu} in Lemma~\ref{le03}, we obtain
\begin{equation}
\label{eqdd}
E(\phi+h)-E(\phi)
\ge k_0\|h\|_{H^1}^2-Ct^2-CY_R(h)-C\|h\|_{H^1}^3-C_R\|h\|_{H^1}^{2+\alpha}.
\end{equation}
Because $0<\alpha<1$ and $\|h\|_{H^1}\le1$,
the cubic term can be absorbed into the final
remainder after increasing $C_R$.
The assumption $E(\phi+h)\le E(\phi)$ and \eqref{eqdd}
implies
\begin{equation}
\label{eqlk}
k_0\|h\|_{H^1}^2
\le C\bigl(t^2+Y_R(h)\bigr)
+C_R\|h\|_{H^1}^{2+\alpha}.
\end{equation}
For each $R\ge1$, choose $\delta_R>0$ below
the fixed thresholds used above and such that
$C_R\delta_R^\alpha\le k_0/2$. Thus, if $\|h\|_{H^1}\le\delta_R$, 
absorbing the last term in \eqref{eqlk} we obtain \eqref{estih},
where $C_*$ is independent of $R$ and of
the choice of $\chi_R$. Finally, since
$z_+(h)+z_-(h)=2\lambda(u,d)_{L^2}$,
we have
\[
(u,d)_{L^2}^2
\le
\frac{|z_+(h)|^2+|z_-(h)|^2}{2\lambda^2}.
\]
Increasing $C_*$ if necessary proves the
second estimate.
\end{proof}

\begin{lemma}[Modulation and projected equations]
\label{le05}
Fix $0<\alpha<2/N$.
There exist $\rho_*>0$ and $C>0$ with the following
property.
Let $I\subset\mathbb R$ be an interval and let $u$ be a
radial solution of \eqref{NLS} belonging to
$C_w(I;H^1(\mathbb R^N))
\cap C(I;L^2(\mathbb R^N))
\cap L^\infty_{\mathrm{loc}}(I;H^1(\mathbb R^N))
\cap W^{1,\infty}_{\mathrm{loc}}(I;H^{-1}(\mathbb R^N))$.
Suppose that, for every $t\in I$, $M(u(t))=M(\phi)$ and
\begin{equation}
\label{cvnb}
\inf_{\gamma\in\mathbb R}
\|u(t)-e^{i\gamma}\phi\|_{H^1}<\rho_*.
\end{equation}
Then there exists a real-valued function $\theta\in W^{1,\infty}_{\mathrm{loc}}(I)$ such that
$u(t)=e^{i(\omega t+\theta(t))}(\phi+h(t))$, where $h=p+iq$,
$h$ has the same regularity as $u$, and
$(q(t),\phi)_{L^2}=0$, $(\phi+p(t),\phi)_{L^2}>0$.
Furthermore,
\begin{equation}
\label{eqjk}
\|h(t)\|_{H^1}
\le C\inf_{\gamma\in\mathbb R}
\|u(t)-e^{i\gamma}\phi\|_{H^1}.
\end{equation}
for every $t\in I$, and
\begin{equation}
\label{estte}
|\dot\theta(t)|\le C\|h(t)\|_{H^1}
\quad\text{for almost every }t\in I.
\end{equation}
The functions $z_\pm(t)=z_\pm(h(t))$ are locally absolutely continuous and satisfy $\dot z_+=\lambda z_++r_+$ and $\dot z_-=-\lambda z_-+r_-$
for almost every $t\in I$, with
\begin{equation}
\label{PrE}
|r_+(t)|+|r_-(t)|
\le C\|h(t)\|_{H^1}^{1+\alpha}.
\end{equation}
\end{lemma}
\begin{proof}
Fix $0<\alpha<2/N$.
We first choose the phase by considering
\[
\mathfrak a(t)=
\int_{\R^N}e^{-i\omega t}u(t,x)\phi(x)\,dx.
\]
From \eqref{cvnb} we infer that $|\mathfrak a(t)|\ge M(\phi)-\rho_*\left\|\phi\right\|_{L^2(\R^N)}$.
Thus, after decreasing $\rho_*$ if necessary,
\begin{equation}
\label{posT}
|\mathfrak a(t)|\ge\frac12M(\phi),
\qquad t\in I.
\end{equation}
Since $\phi\in H^1(\R^N)$ and
$u\in W^{1,\infty}_{\mathrm{loc}}(I;H^{-1}(\R^N))$,
the function $\mathfrak a$ is locally Lipschitz.
It therefore admits a locally Lipschitz real-valued
argument on the interval $I$.
Choose such an argument $\theta$, so that
\[
e^{i\theta(t)}
=\frac{\mathfrak a(t)}{|\mathfrak a(t)|},
\]
and define $h(t)=e^{-i(\omega t+\theta(t))}u(t)-\phi$. Writing
$h=p+iq$, we obtain $(q(t),\phi)_{L^2}=0$ and
$(\phi+p(t),\phi)_{L^2}=|\mathfrak a(t)|>0$.
Multiplication by the locally Lipschitz phase factor
also shows that $h$ has the regularity asserted
in the lemma.

We now prove \eqref{eqjk}. For every $\beta\in\R$, Cauchy--Schwarz gives
\[
|e^{-i\beta}\mathfrak a(t)-M(\phi)|
\le \left\|\phi\right\|_{L^2(\R^N)}
\left\|e^{-i\omega t}u(t)-e^{i\beta}\phi\right\|_{L^2(\R^N)}.
\]
Since $|z/|z|-1|\le2|z-m|/m$ for $z\ne0$ and
$m>0$, the choice
$e^{i\theta(t)}=\mathfrak a(t)/|\mathfrak a(t)|$
implies
\begin{equation}\label{esrb}
|e^{i(\theta(t)-\beta)}-1|
\le C\left\|e^{-i\omega t}u(t)-e^{i\beta}\phi\right\|_{L^2(\R^N)}.
\end{equation}
By \eqref{esrb} and the triangle inequality then gives
$\|h(t)\|_{H^1}
\le
C\|e^{-i\omega t}u(t)-e^{i\beta}\phi\|_{H^1}$.
Taking a sequence of phases approaching the
infimum, we obtain \eqref{eqjk}.
We decrease $\rho_*$ once more so that
$\|h(t)\|_{H^1}\le1$ for every $t\in I$.

Substituting $u=e^{i(\omega t+\theta)}(\phi+h)$
into \eqref{NLS} and using the stationary equation,
we obtain, in $H^{-1}(\R^N)$ for almost every time (cf. Lemma~\ref{le01}),
\begin{equation}
\begin{aligned}
\partial_t p
&=L_-q-\frac{\operatorname{Im}\mathcal N(\phi,h)}{|x|}
+\dot\theta q,\\
\partial_t q
&=-L_+p+\frac{\operatorname{Re}\mathcal N(\phi,h)}{|x|}
-\dot\theta(\phi+p).
\end{aligned}
\label{mo01}
\end{equation}
Differentiating $(q,\phi)_{L^2}=0$ and using
the second equation in \eqref{mo01}, we find
\[
\dot\theta
=
\frac{
-\mathcal F_+(p,\phi)
+\displaystyle\int_{\R^N}
\frac{\operatorname{Re}\mathcal N(\phi,h)\phi}{|x|}\,dx
}{
M(\phi)+(p,\phi)_{L^2}
}
\]
for almost every $t\in I$.
Since $M(\phi)+(p,\phi)_{L^2}=|\mathfrak a(t)| \geq M(\phi)/2$, continuity of the form $\mathcal F_+$ and estimate \eqref{eqn1} give $|\dot\theta|\le C\bigl(\|h\|_{H^1}+\|h\|_{H^1}^{1+\alpha}\bigr)\le C\|h\|_{H^1}$. This proves the phase derivative estimate.

It remains to derive \eqref{PrE}. Pairing the first equation
in \eqref{mo01} with $d$ and using $L_-d=\xi$,
we obtain
\begin{equation}
\label{els}
\lambda\frac{d}{dt}(p,d)_{L^2}
=\lambda(q,\xi)_{L^2}
-\lambda\int_{\mathbb R^N}\frac{\operatorname{Im}\mathcal N(\phi,h)d}{|x|}\,dx
+\lambda\dot\theta(q,d)_{L^2}.
\end{equation}
Similarly, pairing the second equation in \eqref{mo01} with $\xi$ and using $L_+\xi=-\lambda^2d+\kappa\phi$ and $(\phi,\xi)_{L^2}=0$,
gives
\begin{equation}
\label{ecfv}
\frac{d}{dt}(q,\xi)_{L^2}
=\lambda^2(p,d)_{L^2}-\kappa(p,\phi)_{L^2}+\int_{\mathbb R^N}\frac{\operatorname{Re}\mathcal N(\phi,h)\xi}{|x|}\,dx-\dot\theta(p,\xi)_{L^2}.
\end{equation}
The mass assumption implies $M(\phi+h)=M(\phi)$ and $(p,\phi)_{L^2}=-\frac12\left\|h\right\|_{L^2(\R^N)}^2$. Consequently, adding and subtracting the preceding identities \eqref{els}--\eqref{ecfv} yields $\dot z_+=\lambda z_++r_+$ and $\dot z_-=-\lambda z_-+r_-$,
where
\[
\begin{aligned}
r_\pm
={}&
-\lambda\int_{\R^N}
\frac{\operatorname{Im}\mathcal N(\phi,h)d}{|x|}\,dx
+\lambda\dot\theta(q,d)_{L^2}\\
&\quad\pm\left(
\frac{\kappa}{2}\left\|h\right\|_{L^2(\R^N)}^2
+\int_{\R^N}
\frac{\operatorname{Re}\mathcal N(\phi,h)\xi}{|x|}\,dx
-\dot\theta(p,\xi)_{L^2}
\right).
\end{aligned}
\]
By \eqref{eqn1} in Lemma~\ref{le03}, the phase derivative estimate \eqref{estte},
and the Cauchy--Schwarz inequality, we see that $|r_+|+|r_-|
\le
C\bigl(\|h\|_{H^1}^{1+\alpha}+\|h\|_{H^1}^2\bigr)$. Since $0<\alpha<1$ and $\|h\|_{H^1}\le1$, 
we obtain \eqref{PrE}. This completes the proof of the lemma.
\end{proof}

\begin{lemma}[Localized mass identities]\label{le06}
Let $I\subset\R$ be an interval, and let $u$ be
a solution of \eqref{NLS} such that  $u\in C(I;L^2(\R^N))
\cap L^\infty_{\mathrm{loc}}(I;H^1(\R^N))
\cap W^{1,\infty}_{\mathrm{loc}}(I;H^{-1}(\R^N))$.

For every real-valued radial function
$\chi\in C^\infty(\R^N)\cap W^{1,\infty}(\R^N)$,
the function $t\mapsto\int_{\R^N}\chi|u(t)|^2\,dx$
belongs to $W^{1,\infty}_{\mathrm{loc}}(I)$ and
satisfies
\begin{equation}
\label{MassLo}
\frac{d}{dt}\int_{\mathbb{R}^N}\chi |u|^2\,dx
=
2\operatorname{Im}\int_{\mathbb{R}^N}
\nabla\chi\cdot\overline{u}\,\nabla u\,dx.
\end{equation}
for almost every $t\in I$.
Assume in addition that $u$ satisfies the hypotheses
of Lemma~\ref{le05}, and let
$u=e^{i(\omega t+\theta(t))}(\phi+h)$ be the
corresponding modulation.
Fix a radial function $\chi\in C^\infty(\R^N;[0,1])$
such that $\chi=0$ on $B_1$ and $\chi=1$ on
$\R^N\setminus B_2$. For $R\ge1$, set
$\chi_R(x)=\chi(x/R)$ and
\[
Y_R(t):=\int_{\R^N}\chi_R|h(t)|^2\,dx.
\]
Then $Y_R\in W^{1,\infty}_{\mathrm{loc}}(I)$ and
\begin{equation}
\label{YR34}
\begin{split}
Y_R'(t)
={}&2\operatorname{Im}\int_{\mathbb{R}^N}
\nabla\chi_R\cdot\overline h\,\nabla h\,dx
+2\dot\theta\operatorname{Im}\int_{\mathbb{R}^N}
\chi_R\phi\,\overline h\,dx\\
&-2\operatorname{Im}\int_{\mathbb{R}^N}
\frac{\chi_R}{|x|}\bigl(g(\phi+h)-g(\phi)\bigr)\overline h\,dx.
\end{split}
\end{equation}
for almost every $t\in I$. Moreover, writing
$a=1+\sqrt{1+\omega}$, we have
\begin{equation}\label{YRve}
|Y_R'(t)|
\le\varepsilon_R\|h(t)\|_{H^1}^2,
\qquad
\varepsilon_R
=C\left(R^{-1}+R^{(N-1)/2}e^{-aR}\right)
\to0.
\end{equation}
as $R\to\infty$, where $C$ is independent of $R$.
Consequently, for every $s<t$ in $I$,
\begin{equation}
\label{YRSS}
|Y_R(t)-Y_R(s)|
\le \varepsilon_R
\int_s^t\|h(\tau)\|_{H^1}^2\,d\tau.
\end{equation}
\end{lemma}

\begin{proof}
The localized mass identity \eqref{MassLo} is standard; see, for
example, \cite[Section~6]{DHR8}.
At the regularity considered here, it follows by
mollification in time in the $H^{-1}$--$H^1$
duality pairing, as in the proof of
\cite[Lemma~6.5.2]{C03}.
Since $\chi$ is real-valued and
$g(u)\overline u=|u|^2\log|u|^2$ is real,
the nonlinear contribution vanishes, which implies the identity \eqref{MassLo}.

We next establish the identity \eqref{YR34}. Recall that $u=e^{i(\omega t+\theta(t))}(\phi+h)$. Thus, we see that $h$ satisfies
\[
i\partial_t h
=
-\Delta h+\omega h
-\frac{g(\phi+h)-g(\phi)}{|x|}
+\dot\theta(\phi+h)
\]
in $H^{-1}(\R^N)$ for almost every $t\in I$.
Since $Y_R'(t)=2\operatorname{Re}\langle\partial_t h,\chi_Rh\rangle_{H^{-1},H^1}$,
integration by parts yields \eqref{YR34}. We now estimate the three terms in \eqref{YR34}.
Since $\left\|\nabla\chi_R\right\|_{L^\infty(\R^N)}
=R^{-1}\left\|\nabla\chi\right\|_{L^\infty(\R^N)}$, we have
\[
\left|
2\operatorname{Im}\int_{\R^N}
\nabla\chi_R\cdot\overline h\,\nabla h\,dx
\right|
\le
\frac{C}{R}\left\|h\right\|_{L^2(\R^N)}\left\|\nabla h\right\|_{L^2(\R^N)}
\le\frac{C}{R}\|h\|_{H^1}^2.
\]
By \cite[Lemma~9.3.5]{C03},
\[
\left|
\operatorname{Im}
\bigl((g(z)-g(w))\overline{z-w}\bigr)
\right|
\le4|z-w|^2,
\qquad z,w\in\C.
\]
Taking $z=\phi+h$ and $w=\phi$, and using
$\chi_R=0$ on $B_R$, we obtain
\[
\left|
2\operatorname{Im}\int_{\R^N}
\frac{\chi_R}{|x|}
\bigl(g(\phi+h)-g(\phi)\bigr)\overline h\,dx
\right|
\le\frac{8}{R}\left\|h\right\|_{L^2(\R^N)}^2.
\]

Finally, \eqref{estte} in Lemma~\ref{le05} gives
$|\dot\theta|\le C\|h\|_{H^1}$.
Thus, by Cauchy--Schwarz,
\[
\left|
2\dot\theta\operatorname{Im}
\int_{\R^N}\chi_R\phi\,\overline h\,dx
\right|
\le C\left\|\chi_R\phi\right\|_{L^2(\R^N)}\|h\|_{H^1}^2.
\]
The explicit formula
$\phi(x)=e^{(N-1)a/2}e^{-a|x|}$ yields
\[
\left\|\chi_R\phi\right\|_{L^2(\R^N)}^2
\le\int_{|x|>R}\phi^2\,dx
\le C\int_R^\infty r^{N-1}e^{-2ar}\,dr
\le CR^{N-1}e^{-2aR},
\qquad R\ge1.
\]
Combining these estimates proves the bound  for
$Y_R'$. The inequality \eqref{YRSS} follows by integration in time and the fundamental theorem of calculus.
\end{proof}

\begin{proof}[Proof of Theorem~\ref{thm1}]
Fix $-1<\omega<-1+(N-1)^{-2}$ and
$0<\alpha<2/N$. Let $\xi,d,\lambda,\kappa$
be given by Lemma~\ref{sp02}, with
$(\xi,d)_{L^2}=1$.

For radial initial data, the construction in
Proposition~\ref{pr01} yields radial solutions.
Indeed, the regularized equations \eqref{eq07} have radial
coefficients and unique solutions, and their
weak limits remain radial. We first construct initial data with the same
mass as $\phi$ and strictly smaller energy.
Set
\[
v=\xi+i\frac{\lambda}{2}d,
\qquad
s_\varepsilon=
\left(
\frac{M(\phi)}
{M(\phi)+\varepsilon^2\left\|v\right\|_{L^2(\R^N)}^2}
\right)^{1/2},
\qquad
u_{0,\varepsilon}
=s_\varepsilon(\phi+\varepsilon v).
\]
Since $\xi$ and $d$ are orthogonal to $\phi$,
we have $M(u_{0,\varepsilon})=M(\phi)$ and
$u_{0,\varepsilon}\to\phi$ in $H^1(\R^N)$.
Notice also that since $\phi,\xi,d\in H^1_{\mathrm{rad}}(\R^N)$,
we have $u_{0,\varepsilon}\in H^1_{\mathrm{rad}}(\R^N)$.
Writing $h_{0,\varepsilon}=u_{0,\varepsilon}-\phi$, we obtain
$s_\varepsilon=1+O(\varepsilon^2)$ and
$h_{0,\varepsilon}=s_\varepsilon\varepsilon v+(s_\varepsilon-1)\phi$.
In particular, $(\operatorname{Im}h_{0,\varepsilon},\phi)_{L^2}=0$
and $(\phi+\operatorname{Re}h_{0,\varepsilon},\phi)_{L^2}
=s_\varepsilon M(\phi)>0$.

By Lemma~\ref{le02},
$\phi^{-\alpha/(2+\alpha)}v$ and
$\phi^{-\alpha/(2+\alpha)}\phi$ belong to
$H^1(\R^N)$. Consequently, $\left\|\phi^{-\alpha/(2+\alpha)}h_{0,\varepsilon}\right\|_{H^1}\le C\varepsilon$.
Lemma~\ref{le03} and \eqref{eq03} imply
\begin{equation}\label{mate}
\int_{\R^N}\frac{|\mathcal R(\phi,h_{0,\varepsilon})|}{|x|}\,dx
\le C\int_{\R^N}\frac{|\phi^{-\alpha/(2+\alpha)}h_{0,\varepsilon}|^{2+\alpha}}{|x|}\,dx
\le C\varepsilon^{2+\alpha}.
\end{equation}
Moreover, Lemma~\ref{sp02} gives
$\mathcal F_+[\xi]=-\lambda^2$ and
$\mathcal F_-[d]=1$. As $M(u_{0,\varepsilon})=M(\phi)$,
using Lemma~\ref{le07} and \eqref{mate}, we obtain
\begin{equation}\label{in01}
\begin{aligned}
E(u_{0,\varepsilon})-E(\phi)
&=
\frac12\mathcal F_+
[\operatorname{Re}h_{0,\varepsilon}]
+\frac12\mathcal F_-
[\operatorname{Im}h_{0,\varepsilon}]
+O(\varepsilon^{2+\alpha})\\
&=-\frac38\lambda^2\varepsilon^2
+O(\varepsilon^{2+\alpha})<0
\end{aligned}
\end{equation}
for sufficiently small $\varepsilon>0$. 
Let $C_*\ge1$ be the constant in \eqref{vity}, and let
$C_{\mathrm m}\ge1$ and $C_{\mathrm r}>0$ be the constants in
\eqref{eqjk} and \eqref{PrE}, respectively. Thus, by
Lemma~\ref{le05}, we have
\begin{equation}\label{modutes}
\begin{aligned}
&\|h(t)\|_{H^1}
\le C_{\mathrm m}
\inf_{\gamma\in\R}
\|u(t)-e^{i\gamma}\phi\|_{H^1},\\
&|r_+|+|r_-|
\le C_{\mathrm r}\|h\|_{H^1}^{1+\alpha}.
\end{aligned}
\end{equation}
By the definition of $z_\pm$, we have
$z_+(h_{0,\varepsilon})=\frac32\lambda s_\varepsilon\varepsilon$
and $z_-(h_{0,\varepsilon})=\frac12\lambda s_\varepsilon\varepsilon$.
Since $h_{0,\varepsilon}/\varepsilon\to v$
in $L^2(\R^N)$, for each fixed $R$,
\[
\frac{
z_-(h_{0,\varepsilon})^2+Y_R(h_{0,\varepsilon})
}{
z_+(h_{0,\varepsilon})^2
}
\longrightarrow
\frac19+\frac{4}{9\lambda^2}
\int_{\R^N}\chi_R|v|^2\,dx.
\]
The last integral tends to zero as $R\to\infty$.
We may therefore fix $R\ge1$ sufficiently large
that
\begin{equation}\label{in02}
K\varepsilon_R\le\frac{\lambda}{16},
\qquad
\frac19+\frac{4}{9\lambda^2}
\int_{\R^N}\chi_R|v|^2\,dx<\frac14,
\end{equation}
where $\varepsilon_R$ is given in \eqref{YRve} (cf.  Lemma~\ref{le06}) and $K=3C_*/2$.
For this $R$, choose
$0<\delta\le\min\{1,\delta_R\}$ such that
\begin{equation}\label{in03}
C_{\mathrm r}\sqrt K\,\delta^\alpha
\le\frac{\lambda}{16},
\end{equation}
and fix
$0<\rho<\min\{\rho_*,\delta/C_{\mathrm m}\}$,
where $\rho_*$ is given by Lemma~\ref{le05}.
For all sufficiently small $\varepsilon>0$,
\eqref{in01} holds and
\begin{equation}\label{in04}
\|u_{0,\varepsilon}-\phi\|_{H^1}<\rho,
\qquad
z_+(h_{0,\varepsilon})>0,
\qquad
\frac{
z_-(h_{0,\varepsilon})^2+Y_R(h_{0,\varepsilon})
}{
z_+(h_{0,\varepsilon})^2
}<\frac14.
\end{equation}
Fix such an $\varepsilon$, and let $u$ be any
global radial solution with initial datum
$u_{0,\varepsilon}$ satisfying the conclusions
of Proposition~\ref{pr01}.
Suppose, for a contradiction, that
\begin{equation}\label{in05}
\inf_{\gamma\in\R}
\|u(t)-e^{i\gamma}\phi\|_{H^1}<\rho
\qquad\text{for every }t\ge0.
\end{equation}
Lemma~\ref{le05} provides the modulation
$u(t)=e^{i(\omega t+\theta(t))}(\phi+h(t))$
on $[0,\infty)$, with $\theta(0)=0$,
$h(0)=h_{0,\varepsilon}$, and
$\|h(t)\|_{H^1}<\delta$.
Mass conservation, phase invariance, and the energy inequality give
$M(\phi+h(t))=M(\phi)$ and
$E(\phi+h(t))\le E(u_{0,\varepsilon})<E(\phi)$.
Thus Lemmas~\ref{le04}--\ref{le06} imply
\begin{equation}\label{in06}
\begin{gathered}
\|h\|_{H^1}^2
\le C_*\bigl(z_+^2+z_-^2+Y_R\bigr),
\qquad Y_R\ge0,\\
\dot z_+=\lambda z_++r_+,
\qquad
\dot z_-=-\lambda z_-+r_-,\\
|r_+|+|r_-|
\le C_{\mathrm r}\delta^\alpha\|h\|_{H^1},
\qquad
|Y_R'|\le\varepsilon_R\|h\|_{H^1}^2.
\end{gathered}
\end{equation}
The scalar functions $z_+,z_-,Y_R$ are locally
absolutely continuous, and the differential
relations hold almost everywhere.
We now use a continuity argument to show that the
unstable coordinate remains positive and grows
exponentially. Indeed, whenever $z_+(t)>0$, define
\[
W(t):=\frac{z_-^2(t)+Y_R(t)}{z_+^2(t)}.
\]
We claim that
\begin{equation}\label{in07}
z_+(t)\ge z_+(0)e^{\lambda t/2}>0,
\qquad
W(t)\le\frac14,
\qquad t\ge0.
\end{equation}

By \eqref{in04}, $z_+(0)>0$ and $W(0)<1/4$.
Continuity therefore gives an initial interval
on which $z_+>0$ and $W<1/2$.
Let $[0,T_*)$ be the maximal such interval.
On this interval, \eqref{in06} yields $\|h\|_{H^1}^2 \le C_*(1+W)z_+^2 \le Kz_+^2$. Since $\|h\|_{H^1}<\delta$, \eqref{in03} implies $|r_+|+|r_-| \le C_{\mathrm r}\delta^\alpha\|h\|_{H^1} \le C_{\mathrm r}\sqrt K\,\delta^\alpha z_+ \le \frac{\lambda}{16}z_+$.
Consequently,
$\dot z_+=\lambda z_++r_+\ge(\lambda/2)z_+$
almost everywhere on $[0,T_*)$, and integration
gives
\[
z_+(t)\ge z_+(0)e^{\lambda t/2},
\qquad 0\le t<T_*.
\]

We next improve the bound $W<1/2$ to $W\le1/4$.
On every compact subinterval of $[0,T_*)$,
$z_+$ is bounded away from zero, so $W$ is
absolutely continuous. Using the equations
in \eqref{in06}, we obtain almost everywhere
\begin{equation}
\label{Wprime}
\begin{aligned}
W'
&=-2\lambda\frac{z_-^2}{z_+^2}
-2\lambda W
+\frac{2z_-r_-}{z_+^2}
+\frac{Y_R'}{z_+^2}
-2W\frac{r_+}{z_+}\\
&\le
-2\lambda W
+(2\sqrt W+2W)\frac{\lambda}{16}
+K\varepsilon_R.
\end{aligned}
\end{equation}
Here we used $Y_R\ge0$, which implies
$|z_-|/z_+\le\sqrt W$, and
$|Y_R'|\le\varepsilon_R\|h\|_{H^1}^2
\le K\varepsilon_Rz_+^2$.
By \eqref{in02} and \eqref{Wprime}, whenever $1/4\le W<1/2$,
\[
W'\le-\frac{\lambda}{2}
+\frac{\lambda}{16}(\sqrt2+2)
\le-\frac{\lambda}{4}.
\]
Thus $(W-\tfrac14)_+ = \max\left\{W-\frac14,\,0\right\}$ is absolutely continuous on
each compact subinterval of $[0,T_*)$ and
has nonpositive derivative almost everywhere.
Since it vanishes at $t=0$, it remains zero.
Hence $W(t)\le\frac14$ for all $0\le t<T_*$.

If $T_*<\infty$, continuity gives (cf. Lemma \ref{le05})
\[
z_+(T_*)\ge z_+(0)e^{\lambda T_*/2}>0,
\qquad
W(T_*)\le\frac14<\frac12.
\]
These strict inequalities ensure that $z_+>0$
and $W<1/2$ also hold on an interval beyond
$T_*$, contradicting maximality.
Therefore $T_*=\infty$, which proves
\eqref{in07}.

On the other hand, \eqref{in05} and the definition
of $z_+$ imply
\[
|z_+(t)|
\le
\bigl(\lambda\left\|d\right\|_{L^2(\R^N)}+\left\|\xi\right\|_{L^2(\R^N)}\bigr)
\|h(t)\|_{H^1}
\le
\bigl(\lambda\left\|d\right\|_{L^2(\R^N)}+\left\|\xi\right\|_{L^2(\R^N)}\bigr)\delta
\qquad\text{for every }t\ge0.
\]
This contradicts \eqref{in07} as $t\to\infty$,
since $z_+(0)>0$.
Therefore, there exists $t_\varepsilon>0$ such that
\[
\inf_{\gamma\in\R}
\|u(t_\varepsilon)-e^{i\gamma}\phi\|_{H^1}
\ge\rho.
\]
Since $u_{0,\varepsilon}\to\phi$ in $H^1(\R^N)$
and $\rho>0$ is independent of $\varepsilon$,
the standing wave $e^{i\omega t}\phi$ is orbitally
unstable. This proves the theorem.
\end{proof}

\appendix
\section{Spectral properties of the linearized operator}\label{asp}

Throughout this appendix, we assume $N\ge2$ and fix
$-1<\omega<-1+{(N-1)^{-2}}$,
$a:=1+\sqrt{1+\omega}$,
and write $\phi(x)=e^{(N-1)a/2}e^{-a|x|}$.
All operators and quadratic forms below are considered
on the real radial space $L^2_{\mathrm{rad}}(\mathbb{R}^N)$.
The estimate \eqref{ne02} and Young's inequality give
\begin{equation}\label{Hty}
\int_{\R^N}\frac{|v|^2}{|x|}\,dx
\le \frac{2}{N-1}\left\|v\right\|_{L^2(\R^N)}\left\|\nabla v\right\|_{L^2(\R^N)}
\le \varepsilon\left\|\nabla v\right\|_{L^2(\R^N)}^2
+C_\varepsilon\left\|v\right\|_{L^2(\R^N)}^2,
\qquad v\in H^1(\R^N).
\end{equation}
for every $\varepsilon>0$.
Consequently, the quadratic forms
\[
\begin{aligned}
\mathcal F_-[v]
&:=\left\|\nabla v\right\|_{L^2(\R^N)}^2+a^2\left\|v\right\|_{L^2(\R^N)}^2
-(N-1)a\int_{\R^N}\frac{|v|^2}{|x|}\,dx,\\
\mathcal F_+[v]
&:=\left\|\nabla v\right\|_{L^2(\R^N)}^2+a^2\left\|v\right\|_{L^2(\R^N)}^2
-\bigl((N-1)a+2\bigr)
\int_{\R^N}\frac{|v|^2}{|x|}\,dx
\end{aligned}
\]
are closed and bounded below on the common form domain
$H^1_{\mathrm{rad}}(\R^N)$.
We denote their associated self-adjoint operators by
\[
L_-=-\Delta+a^2-\frac{(N-1)a}{|x|},
\qquad
L_+=-\Delta+a^2-\frac{(N-1)a+2}{|x|}.
\]
The compactness established in Lemma~\ref{CCD},
together with interpolation, implies that
$v\mapsto |x|^{-1/2}v$ is compact from
$H^1(\R^N)$ into $L^2(\R^N)$.
Thus the Coulomb terms are relatively compact
perturbations in the form sense, and
$ \sigma_{\mathrm{ess}}(L_-)
=\sigma_{\mathrm{ess}}(L_+)
=[a^2,\infty)$.

\begin{lemma}\label{sp01}
Let $N\ge2$ and $-1<\omega<-1+(N-1)^{-2}$.
Then the following assertions hold.
\begin{enumerate}[label=\textup{(\roman*)}]
\item
The operator $L_-$ is nonnegative and
$\ker L_-=\operatorname{span}\{\phi\}$.
Moreover, there exists $c>0$ such that, for every
$v\in H^1_{\mathrm{rad}}(\R^N)$ satisfying
$(v,\phi)_{L^2}=0$, we have
\begin{equation}\label{cm01}
\mathcal F_-[v]\ge c\|v\|_{H^1}^2.
\end{equation}

\item
The operator $L_+$ has exactly one negative
eigenvalue, which is simple, and $\ker L_+=\{0\}$.
\end{enumerate}
\end{lemma}

\begin{proof}
Recall that $a=1+\sqrt{1+\omega}$, so that $1<a<N/(N-1)$.
We first recall the radial Coulomb eigenvalues.
For $\gamma>0$, the eigenvalues of
$-\Delta-\gamma/|x|$ below zero, with the
$H^1_{\mathrm{rad}}$ form realization, are (cf. \cite[Section~II, equations~(11) and~(16)]{KLPW})
\[
-\frac{\gamma^2}{(2j+N-1)^2},
\qquad j=0,1,2,\ldots,
\]
and each is simple. Applying this formula with $\gamma=(N-1)a$,
the eigenvalues of $L_-$ below $a^2$ are
\[
a^2-\frac{(N-1)^2a^2}{(2j+N-1)^2},
\qquad j=0,1,2,\ldots.
\]
The lowest eigenvalue is zero and is simple,
while all the others are positive.
Since $L_-\phi=0$ and
$\sigma_{\mathrm{ess}}(L_-)=[a^2,\infty)$,
we obtain
$L_-\ge0$ and
$\ker L_-=\operatorname{span}\{\phi\}$.

Consequently, $\mathcal F_-[v]>0$ for every
nonzero $v\in H^1_{\mathrm{rad}}(\R^N)$
orthogonal to $\phi$.
Indeed, if $\mathcal F_-[v]=0$, nonnegativity
of the form implies $v\in\ker L_-$, and the
orthogonality condition then gives $v=0$.
Thus, Lemma \ref{CCD} and a standard compactness argument yield \textup{(i)}.
Finally, applying the radial Coulomb formula
with $\gamma=(N-1)a+2$, the eigenvalues of $L_+$
below $a^2$ are
\[
a^2-\frac{\bigl((N-1)a+2\bigr)^2}
{(2j+N-1)^2},
\qquad j=0,1,2,\ldots.
\]
For $j=0$, this eigenvalue is strictly negative.
For $j=1$, it is strictly positive because $(N-1)a+2<(N+1)a$,
which follows from $a>1$.
The remaining eigenvalues are larger, and the
essential spectrum starts at $a^2>0$.
Thus $L_+$ has exactly one negative eigenvalue,
which is simple, and zero does not belong to
its spectrum. This proves \textup{(ii)}.
\end{proof}

For the following result, we use the notation
$X:=\{v\in L^2_{\mathrm{rad}}(\R^N): (v,\phi)_{L^2}=0\}$.
\begin{lemma}\label{sp02}
Let $N\ge2$ and $-1<\omega<-1+(N-1)^{-2}$.
The linearized operator
\[
\mathcal A=
\begin{pmatrix}
0&L_-\\
-L_+&0
\end{pmatrix},
\qquad
D(\mathcal A)=D(L_+)\times D(L_-),
\]
has a pair of real eigenvalues $\pm\lambda$,
with $\lambda>0$.
These are its only eigenvalues outside the
imaginary axis in the complexified radial space.

There exist real-valued functions
$\xi\in D(L_+)\cap X$ and $d\in D(L_-)\cap X$,
and a constant $\kappa\in\R$, such that
\begin{equation}\label{df}
L_-d=\xi,\qquad
L_+\xi=-\lambda^2d+\kappa\phi,\qquad
(\xi,d)_{L^2}=1.
\end{equation}
In particular, setting
$\eta=\lambda d-\kappa\phi/\lambda$, we have
$\eta\in D(L_-)$ and
\[
\mathcal A
\begin{pmatrix}\xi\\ \pm\eta\end{pmatrix}
=\pm\lambda
\begin{pmatrix}\xi\\ \pm\eta\end{pmatrix}.
\]
\end{lemma}

\begin{proof}
By Lemma~\ref{sp01}, $L_-$ is nonnegative,
$\ker L_-=\operatorname{span}\{\phi\}$, and
its form is coercive on
$H^1_{\mathrm{rad}}(\R^N)\cap X$.
Thus $L_-^{1/2}$ is well defined, with $D(L_-^{1/2})=H^1_{\mathrm{rad}}(\R^N)$.

Moreover, the spectrum of $L_-$ on $X$ is
bounded away from zero (cf. Lemma \ref{sp01}(i)).
We shall use $L_-^{-1}$ and $L_-^{-1/2}$
only on $X$, where they are bounded operators
defined by the spectral calculus.
Consider the quadratic form
\[
\mathfrak k[f]
:=
\mathcal F_+[L_-^{1/2}f]
=
\left\|L_-f\right\|_{L^2(\R^N)}^2
-2\int_{\R^N}
\frac{|L_-^{1/2}f|^2}{|x|}\,dx,
\qquad f\in D(L_-).
\]
We denote its associated self-adjoint operator
by $\mathcal K$.
Thus $\mathcal K$ is the form realization of
$L_-^{1/2}L_+L_-^{1/2}$.
To justify this construction, the
estimate \eqref{Hty} and equivalence of the shifted form
norm of $L_-$ with the $H^1$ norm give
\[
\int_{\R^N}
\frac{|L_-^{1/2}f|^2}{|x|}\,dx
\le
\varepsilon\left\|L_-f\right\|_{L^2(\R^N)}^2
+C_\varepsilon\left\|L_-^{1/2}f\right\|_{L^2(\R^N)}^2.
\]
Spectral interpolation yields $\left\|L_-^{1/2}f\right\|_{L^2(\R^N)}^2
\le
\delta\left\|L_-f\right\|_{L^2(\R^N)}^2+C_\delta\left\|f\right\|_{L^2(\R^N)}^2.$
Consequently, the Coulomb term is infinitesimally
form bounded relative to $\left\|L_-f\right\|_{L^2(\R^N)}^2$.
Hence $\mathfrak k$ is closed and bounded below.

Furthermore, $L_-^{1/2}$ maps $D(L_-)$,
equipped with its graph norm, continuously
into $H^1_{\mathrm{rad}}(\R^N)$.
The compactness of
$v\mapsto |x|^{-1/2}v$ from $H^1$ into $L^2$ (cf. Lemma \ref{CCD})
therefore shows that $\mathfrak k$ is a compact
form perturbation of the form of $L_-^2$.
It follows that $ \sigma_{\mathrm{ess}}(\mathcal K)
=\sigma_{\mathrm{ess}}(L_-^2)
=[a^4,\infty).$

We next show that $\mathcal K$ has a negative
eigenvalue. Set $\zeta(x)=\left(\tfrac N2-a|x|\right)\phi(x).$
Direct integration gives
$(\zeta,\phi)_{L^2}=0$, so that
$\zeta\in H^1_{\mathrm{rad}}(\R^N)\cap X$.
Using the explicit exponential profile, we obtain
\[
\mathcal F_-[\zeta]=a^2M(\phi),
\qquad
\int_{\R^N}\frac{|\zeta|^2}{|x|}\,dx
=\frac N4m_1(\phi),
\qquad
m_1(\phi)=\frac{2a}{N-1}M(\phi).
\]
Therefore,
\[
\mathcal F_+[\zeta]
=a^2M(\phi)-\frac N2m_1(\phi)
=aM(\phi)\left(a-\frac{N}{N-1}\right)<0.
\]
Since $\zeta\in D(L_-^{1/2})\cap X$,
the spectral calculus gives
$f_0:=L_-^{-1/2}\zeta\in D(L_-)$, and $\mathfrak k[f_0]=\mathcal F_+[\zeta]<0$.
As the essential spectrum of $\mathcal K$
is contained in $[a^4,\infty)$, this implies the existence of a negative
eigenvalue.

Lemma~\ref{sp01} also shows that $\mathcal F_+$
has exactly one negative direction.
Since
$\mathfrak k[f]=\mathcal F_+[L_-^{1/2}f]$,
the form $\mathfrak k$ has at most one negative
direction: on any negative subspace, the map
$f\mapsto L_-^{1/2}f$ is injective and preserves
negativity.
Thus $\mathcal K$ has exactly one negative
eigenvalue, which is simple.
Write it as $-\lambda^2$, with $\lambda>0$,
and choose a nonzero real eigenfunction
$f\in D(\mathcal K)$.

Note that $\phi$ belongs to $\ker\mathcal K$,
because $L_-^{1/2}\phi=0$.
Self-adjointness therefore gives $(f,\phi)_{L^2}=0$.
Set
\[
\xi=L_-^{1/2}f,\qquad d=L_-^{-1}\xi.
\]
Then $\xi\in H^1_{\mathrm{rad}}(\R^N)\cap X$,
$d\in D(L_-)\cap X$, and $L_-d=\xi$.
For every $v\in H^1_{\mathrm{rad}}(\R^N)\cap X$,
the function $L_-^{-1/2}v$ belongs to $D(L_-)$.
Testing the form eigenvalue equation against
this function yields
\[
\mathcal F_+(\xi,v)
=-\lambda^2(f,L_-^{-1/2}v)_{L^2}
=-\lambda^2(d,v)_{L^2}.
\]
Here $\mathcal F_+(\cdot,\cdot)$ denotes the
symmetric bilinear form associated with
$\mathcal F_+$.
Since $H^1_{\mathrm{rad}}(\R^N)\cap X$ has
codimension one in $H^1_{\mathrm{rad}}(\R^N)$,
there exists $\kappa\in\R$ such that $L_+\xi=-\lambda^2d+\kappa\phi$
in the form sense.
The right-hand side belongs to $L^2$, so
$\xi\in D(L_+)$ and the identity also holds
in the operator sense.

Moreover, $(\xi,d)_{L^2}
=(\xi,L_-^{-1}\xi)_{L^2}>0.$
Rescaling $\xi$ and $d$ by the same positive
constant, and adjusting $\kappa$ accordingly,
we may impose $(\xi,d)_{L^2}=1$.
Define $\eta=\lambda d-\frac{\kappa}{\lambda}\phi.$
Then $\eta\in D(L_-)$ and
\[
L_-\eta=\lambda\xi,\qquad
-L_+\xi=\lambda\eta.
\]
Thus $(\xi,\eta)$ is an eigenvector of
$\mathcal A$ with eigenvalue $\lambda$.
Replacing $\eta$ by $-\eta$ gives the
eigenvalue $-\lambda$.
Finally, let $z\ne0$ be an eigenvalue of
$\mathcal A$, with eigenvector $(p,q)$ in
the complexified domain.
Then
\[
L_-q=zp,\qquad -L_+p=zq.
\]
Pairing the first equation with $\phi$ gives
$(p,\phi)_{L^2}=0$.
Also, $p\ne0$, since otherwise the second
equation would give $q=0$.
The first equation implies
$q=zL_-^{-1}p+c\phi$
for some $c\in\C$.
Substituting this identity into the second
equation gives
$L_+p=-z^2L_-^{-1}p-zc\phi$.
Set $F=L_-^{-1/2}p\in D(L_-)$.
For every $g\in D(L_-)$, the vector
$L_-^{1/2}g$ is orthogonal to $\phi$.
Using the sesquilinear extensions of the forms,
we obtain
\[
\mathfrak k(F,g)
=\mathcal F_+(p,L_-^{1/2}g)
=-z^2(F,g)_{L^2}.
\] Therefore $F\in D(\mathcal K)$ and
$\mathcal K F=-z^2F$.
Self-adjointness of $\mathcal K$ implies
$z^2\in\R$, so $z$ is either real or purely
imaginary.
If $z$ lies outside the imaginary axis,
then $-z^2<0$.
Since $-\lambda^2$ is the only negative
eigenvalue of $\mathcal K$, we conclude
that $z=\pm\lambda$.
\end{proof}

\begin{lemma}\label{sp03}
Let $N\ge2$ and $-1<\omega<-1+(N-1)^{-2}$,
and let $d$ (cf. \eqref{df}) be the function given by Lemma~\ref{sp02}.
There exists $c>0$ such that, for every
$v\in H^1_{\mathrm{rad}}(\R^N)$ satisfying $(v,\phi)_{L^2}=(v,d)_{L^2}=0$,
we have
\begin{equation}\label{cp01}
\mathcal F_+[v]\ge c\|v\|_{H^1}^2.
\end{equation}
\end{lemma}

\begin{proof}
Set
\[
V:=\{v\in H^1_{\mathrm{rad}}(\R^N):
(v,\phi)_{L^2}=(v,d)_{L^2}=0\}.
\]
We first prove that $\mathcal F_+[v]>0$
for every nonzero $v\in V$.
By Lemma~\ref{sp02},
\[
L_+\xi=-\lambda^2d+\kappa\phi,
\qquad
(\xi,\phi)_{L^2}=0,
\qquad
(\xi,d)_{L^2}=1.
\]
Consequently,
\[
\mathcal F_+[\xi]=-\lambda^2,
\qquad
\mathcal F_+(\xi,v)=0
\quad\text{for every }v\in V.
\]
If $\mathcal F_+[v]<0$ for some $v\in V$,
then $\xi$ and $v$ are linearly independent
and $\mathcal F_+$ is negative definite on
$\operatorname{span}\{\xi,v\}$.
This contradicts Lemma~\ref{sp01}, which states
that $L_+$ has exactly one negative eigenvalue,
counted with multiplicity.
Thus $\mathcal F_+$ is nonnegative on $V$.
Suppose now that $v\in V$ and $\mathcal F_+[v]=0$.
Nonnegativity of the restricted form implies
$\mathcal F_+(v,w)=0$ for every $w\in V$.
The same pairing vanishes for $w=\xi$.
Since every $w\in H^1_{\mathrm{rad}}(\R^N)\cap X$
can be written as
\[
w=(w,d)_{L^2}\xi+
\bigl(w-(w,d)_{L^2}\xi\bigr),
\]
with the second term in $V$, we obtain
$\mathcal F_+(v,w)=0$ for every such $w$.
Hence, for some $\gamma\in\R$, $L_+v=\gamma\phi$.
Since the right-hand side belongs to $L^2$,
we have $v\in D(L_+)$.

By Lemma~\ref{sp01} and
$\sigma_{\mathrm{ess}}(L_+)=[a^2,\infty)$,
zero belongs to the resolvent set of $L_+$.
Differentiating the stationary equation with
respect to $\omega$ gives $L_+\partial_\omega\phi=-\phi.$
This identity holds in the form sense and,
since $\partial_\omega\phi\in H^1_{\mathrm{rad}}$
and $\phi\in L^2$, also in the operator sense.
Therefore,
\[
(L_+^{-1}\phi,\phi)_{L^2}
=-\frac12\frac{d}{d\omega}M(\phi).
\]
The explicit profile yields
\[
M(\phi)
=\frac{|\mathbb S^{N-1}|\Gamma(N)}{(2a)^N}
e^{(N-1)a},
\qquad
\frac{da}{d\omega}=\frac{1}{2(a-1)},
\]
and consequently
\[
(L_+^{-1}\phi,\phi)_{L^2}
=\frac{M(\phi)}{4a(a-1)}
\bigl(N-(N-1)a\bigr)>0,
\]
because $1<a<N/(N-1)$.
Since $v=\gamma L_+^{-1}\phi$ and $(v,\phi)_{L^2}=0$,
we conclude that $\gamma=0$ and hence $v=0$.
This proves strict positivity on $V\setminus\{0\}$. Thus, Lemma \ref{CCD} and a standard compactness argument yield \eqref{cp01}.
\end{proof}

We fix the functions $\xi,d,\eta$ and the constants
$\lambda,\kappa$ supplied by Lemma~\ref{sp02},
with the normalization $(\xi,d)_{L^2}=1$.

\begin{lemma}[Exponential decay of the spectral modes]
\label{le02}
For every $0<\sigma<a$, we have
$e^{\sigma|x|}\xi,\ e^{\sigma|x|}\eta,\ e^{\sigma|x|}d
\in H^1(\R^N)$.
Consequently, for every $0<\alpha<1$,
$\phi^{-\alpha}\psi\in H^1(\R^N)$
for $\psi\in\{\phi,\xi,d\}$.
\end{lemma}

\begin{proof}
For $N\ge 3$, the result follows from \cite[Theorem~1.5]{HL07}. We give a direct proof that also covers $N=2$.
Fix $0<\sigma<a$.
By Lemma~\ref{sp02}, the real-valued functions
$\xi$ and $\eta$ satisfy
\begin{equation*}
-\Delta\xi+a^2\xi=\frac{(N-1)a+2}{|x|}\xi-\lambda\eta,\qquad
-\Delta\eta+a^2\eta=\frac{(N-1)a}{|x|}\eta+\lambda\xi.
\end{equation*}
in the weak sense. Let $R>1$, to be chosen below, and choose a smooth
radial function $\chi$ such that $0\le\chi\le1$,
$\chi=0$ on $\{|x|\le R\}$, and $\chi=1$ on
$\{|x|\ge R+1\}$.
For $L>R+1$, define
\[
f_L(x)=\chi(x)e^{\sigma\min\{|x|,L\}}.
\]
For each fixed $L$, the function $f_L$ is bounded
and Lipschitz. Hence $f_L^2\eta$, $f_L^2\xi\in H^1$.

Testing the two equations and adding, the terms
containing $\lambda$ cancel. Thus,
we obtain
\begin{equation}
\label{edc}
\begin{aligned}
&\left\|\nabla(f_L\xi)\right\|_{L^2(\R^N)}^2+\left\|\nabla(f_L\eta)\right\|_{L^2(\R^N)}^2
+a^2\int_{\R^N}f_L^2(\xi^2+\eta^2)\,dx\\
&=\int_{\R^N}|\nabla f_L|^2(\xi^2+\eta^2)\,dx
+\int_{\R^N}\frac{f_L^2}{|x|}
\left(\bigl((N-1)a+2\bigr)\xi^2+(N-1)a\eta^2\right)\,dx.
\end{aligned}
\end{equation}
For every $\varepsilon>0$, $|\nabla f_L|^2
\le
(1+\varepsilon)\sigma^2f_L^2
+C_{\varepsilon,R}\mathbf 1_{\{R<|x|<R+1\}}$,
where $C_{\varepsilon,R}$ is independent of $L$.
Since $f_L$ vanishes on $\{|x|\le R\}$, by \eqref{edc} it follows
that
\begin{equation*}
\begin{aligned}
&\left\|\nabla(f_L\xi)\right\|_{L^2(\R^N)}^2+\left\|\nabla(f_L\eta)\right\|_{L^2(\R^N)}^2
+\left(a^2-(1+\varepsilon)\sigma^2-\frac{(N-1)a+2}{R}\right)
\int_{\R^N}f_L^2(\xi^2+\eta^2)\,dx\\
&\qquad\le C_{\varepsilon,R}
\int_{R<|x|<R+1}(\xi^2+\eta^2)\,dx.
\end{aligned}
\end{equation*}
Because $\sigma<a$, we may first choose
$\varepsilon>0$ sufficiently small and then
$R$ sufficiently large so that the coefficient
in parentheses is positive.
Consequently,
\begin{equation}\label{eq123}
\sup_{L>R+1}
\left(\|f_L\xi\|_{H^1}+\|f_L\eta\|_{H^1}\right)
<\infty.
\end{equation}
By \eqref{eq123}, weak compactness in $H^1$, together with local
convergence of the truncated weights, yields
$\chi e^{\sigma|x|}\xi,\chi e^{\sigma|x|}\eta
\in H^1(\R^N)$.
Since $(1-\chi)e^{\sigma|x|}$ is bounded and
Lipschitz with compact support, we conclude that
$e^{\sigma|x|}\xi,\ e^{\sigma|x|}\eta\in H^1(\R^N)$.
The explicit formula
$\phi(x)=e^{(N-1)a/2}e^{-a|x|}$ also gives
$e^{\sigma|x|}\phi\in H^1(\R^N)$.
Moreover, Lemma~\ref{sp02} gives $d=\frac{\eta}{\lambda}
+\frac{\kappa}{\lambda^2}\phi$. 
Hence $e^{\sigma|x|}d\in H^1(\R^N)$,
which proves the first assertion.
Finally, let $0<\alpha<1$ and choose
$\alpha a<\sigma<a$.
For $\psi\in\{\phi,\xi,d\}$, we have
\[
\phi^{-\alpha}\psi
=
e^{-\alpha(N-1)a/2}
e^{-(\sigma-\alpha a)|x|}
\bigl(e^{\sigma|x|}\psi\bigr).
\]
The last factor belongs to $H^1$, and the
preceding spatial factor is bounded and
Lipschitz. Thus $\phi^{-\alpha}\psi\in H^1(\R^N)$,
as claimed.
\end{proof}



\end{document}